\documentclass[hidelinks,onefignum,onetabnum]{siamart220329}

\usepackage{lipsum}
\usepackage{amsfonts}
\usepackage{graphicx}
\usepackage{epstopdf}
\usepackage{algorithmic}
\ifpdf
  \DeclareGraphicsExtensions{.eps,.pdf,.png,.jpg}
\else
  \DeclareGraphicsExtensions{.eps}
\fi

\newsiamremark{remark}{Remark}
\newsiamremark{hypothesis}{Hypothesis}
\crefname{hypothesis}{Hypothesis}{Hypotheses}
\newsiamthm{claim}{Claim}

\headers{Pipline flow of blended gas}{Martin Gugat}

\title{Synchronization of velocities in pipline flow of blended gas
\thanks{.
\funding{
This work was funded by DFG in the 
 Collaborative Research Centre
CRC/Transregio 154,
Mathematical Modelling, Simulation and Optimization Using the Example of Gas Networks,
Project C03 and   C05,    Projektnummer 239904186.
}}}

\author{Martin Gugat\thanks{
Friedrich-Alexander-Universit\"at Erlangen-N\"urnberg (FAU),
Chair for Dynamics, Control, Machine Learning and Numerics 
– Alexander von Humboldt Professorship, 
\\
Department Mathematik, Cauerstr. 11, 91058 Erlangen, Germany, 
(\email{martin.gugat@fau.de})
}
}

\usepackage{amsopn}

\ifpdf
\hypersetup{
  pdftitle={synchronization
   MANUSCRIPT EDITION !!!!!!
  },
  pdfauthor={M. Gugat}
}
\fi

\usepackage{float}

\usepackage{yfonts}

\usepackage{tikz} 
\usetikzlibrary{positioning}

\usepackage{amsfonts,amsmath,amssymb}
\usepackage{graphicx}

\renewcommand{\geq}{\geqslant}
\renewcommand{\leq}{\leqslant}

\usepackage{siunitx}

\newtheorem{example}{Example}

\usepackage{color}

\def\rev#1{{\color{black} #1}}

\newcommand{\mar}{\color{black}}

\newcommand{\tMz}{\mbox{$\hat M$} }

\begin{document}

\maketitle

%\begin{abstract}
%\textbf{Abstract:}

%%%%%%%%%%%%%%%%%%%%%%%%%%%%%%%%%%%%%%%%%%%%%%%%%%%%%%%%%%%%%%%%%%%%%%%
%Pipeline networks  can be modelled by  graphs where the edges are
%given by the pipes that form the network.
%The flow of gas through the pipes is governed by 
%a quasi-linear hyperbolic system of partial differential equations,
%namely  the isothermal Euler equations.
%%hyperbolic systems of partial differential equations (PDEs). 
%An additional conservation law
%allows to obtain a model for
%the flow with hydrogen blending.
%Algebraic node conditions define 
%a model for the flow through the 
%vertices  where the pipes are connected.
%%network is obtained 
%%by coupling the PDEs that describe the flow through
%%the 
%%%the isothermal Euler equations
%%by algebraic node conditions
%%that model the flow through the vertices of the  graph. 
%At certain points in space measurements of the state in the network 
%are available. These measurements include a certain measurement error.
%%Based on these nodal observations,
%%the full system state can be approximated using an observer system. 
%In this paper we present
%an observer system for general graphs
%and prove that 
%a nodal coupling
%to the original system 
%drives the 
%state of the observer system
%%converges 
%to the original state
%exponentially fast until
%an error level is reached that 
%depends on the measurement error.
%% Numerical experiments confirm the theoretical findings. 
%\end{abstract}

%%%%%%%%%%%%%%%%%%%%%%%%%%%%%%%%%%%%%%%%%%%%%%%%%%%%%%%%%%%%%%%%%%%%%%%%%%%%%%%%%%%%%%%%%%%%%%%%%%%%%%%%%%%%%%%%%%%%%%%%%%%%%%%%%%%%%%%%%%%%%%%%%%%

\begin{abstract}
% MANUSCRIPT EDITION !!!!!!
We consider the pipeline flow of blended gas.
%, of two gas  components.
The flow  is governed by a coupled system where for each component we have
the isothermal Euler equations with an additional velocity coupling term
that couples the velocities of the different components.  
Our motivation is hydrogen blending in natural gas pipelines, which will play a role in the transition to renewable energies.
We show that with suitable boundary conditions 
%for which  %We show that 
the velocities of the  gas components synchronize exponentially fast,
as long as the $L^2$-norm of the synchronization error is
outside of a certain 
interval where the size of the interval
is determined by the order of the interaction terms.
%until the $L^2$-norm of the difference is on a level that is of the order of the interaction terms.
%Of course this requires that no additional velocity difference is
%introduced through the boundary.
This indicates 
%which shows 
that in some cases  for a mixture of $n$ components it is 
justified 
%sufficient 
to use a 
%simplified $
%(n+1) \times (n+1) $ 
drift-flux 
model where it is assumed that all components flow with the same velocity.
%
%The networked system is defined  on a graph, where 
%the flow through the  vertices  where the pipes are connected
%is governed by algebraic node conditions
%that require the conservation of mass and
%the continuity of the specific enthalpy.
%
For the proofs we use an appropriately chosen Lyapunov function which is based upon the idea of relative energy.

% MANUSCRIPT EDITION !!!!!!
%We consider blended gas pipeline flow.
%, of two gas  components.
%It is governed by a coupled system where for each component we have
%the isothermal Euler equation with an additional velocity coupling term.
%Our motivation is hydrogen blending in natural gas pipelines, which will play a role in the transition to renewable energies.
%We show that the velocities of the two gas components synchronize exponentially fast 
%until the $L^2$-norm of the difference is on a level that is of the order of the interaction terms.
%%Of course this requires that no additional velocity difference is
%introduced through the boundary.
%This indicates 
%which shows 
%that often for a mixture of $n$ components it is sufficient to use a simplified $
%(n+1) \times (n+1) $ model,
%where it is assumed that all components flow with the same velocity.
%
%The networked system is defined  on a graph, where 
%the flow through the  vertices  where the pipes are connected
%is governed by algebraic node conditions
%that require the conservation of mass and
%the continuity of the specific enthalpy.
%
%For the proofs we use an appropriately chosen Lyapunov function which is based upon the idea of relative energy.
%(see e.g. \cite{dafermosbook},  \cite{diperna})
%We show that asymptotically the  velocities of the components synchronize exponentially fast 
%up to the order of the interaction terms. 
%The proof of the synchronization result uses
%Lyapunov functions with exponential weights.
% Numerical experiments confirm the theoretical findings. 
\end{abstract}


%MSC[2010] 

%\end{frontmatter}

%\pagestyle{myheadings}
%\thispagestyle{plain}

%%%%%%%%%%%%%%%%%%%%%%%%%%%%%%%%%%%%%%%%%%%%%%%%%%%%%%%%%%%%%%%%%%%%%%%%%%%

%\maketitle

%\begin{AMS}
%35L04, 49K20, 90C46
%\end{AMS}

%\\
%This work was
%supported by DFG
%in the framework of the Collaborative Research Centre
%CRC/Transregio 154,
%Mathematical Modelling, Simulation and Optimization Using the Example of Gas Networks,
%Project C03 and  Project  C05.

%%%%%%%%%%%%%%%%%%%%%%%%%%%%%%%%%%%%%%%%%%%%%%%%%%%%%%%%%%%%%%%%%%%%%%%%%%%%%%%%%
%%-----------------------------
\section*{Introduction}
The European Union
has  put forward a comprehensive framework to support the uptake of renewable and low-carbon hydrogen to help decarbonise the EU. 
%In the hydrogen policy framework one of the key actions is
Blending hydrogen with natural  gas in pipelines is a part of this strategy.

%%%%%%%%%%%%%%%%%%%%%%%%%%%%%%%%%%%%%%%%%%%%%%%%%%%%%%%%%%%%%%%%%%%%%%%%%%%%%%%%%%%%%%%%%%

As a model for the the flow of the gas mixture
we consider a model where the
flow of each component is governed by the isothermal Euler equations,
with an additional coupling term that models the interactions between
the different components.
For a flow with $n$ components this yields a $2n\times 2n$ quasilinear hyperbolic system
for which Riemann invariants do exist.

We present boundary conditions such that the velocities of the different components synchronize exponentially 
fast with respect to the $L^2$-norm
%We show that for a more general class of boundary conditions 
%the velocities of the  components synchronize exponentially fast
%with respect to the $L^2$-norm 
until the synchronization error
is of an order that is determined by the size of the interaction 
term. 
The proof of the exponential synchronization 
is based on  a suitably chosen Lyapunov function.
Our choice of the Lyapunov function is
motivated by relative energy techniques that have
been used for example in \cite{Dafermos_book}, \cite{herbertegger}.
However, in contrast to \cite{Dafermos_book}, \cite{herbertegger}, 
in our Lyapunov function the derivative of the energy functional does not appear.
The choice of the Lyapunov function is motivated by the fact the
we investigate only the synchronization of the velocities of
the different components, and not of the full state.
For synchronization of the full state,
Lyapunov functions with exponential weights have been used
successfully,
see
\cite{cobada}, \cite{gugiob},
\cite{gugat2011gas}.

The synchronization result indicates that in some cases it  makes 
sense to use a simplified 
 $(n+1)\times (n+1) $ system 
of balance laws as a model,
where it is assumed from
the beginning that the velocities of the components are equal.
For the case of two components, 
such a drift-flux model has been studied 
for example 
in
\cite{driftfluxbanda},
\cite{gugiob}.

%
%In order to point out the difficulties of the analysis, we recap the situation from \cite{castillo2013boundary}.
%Let $\Lambda$ denote a continuously differentiable diagonal matrix function to $[0, \infty)^{n\times n}$
%that contains the positive eigenvalues of the system as diagonal elements.
%A quasi-linear hyperbolic system in diagonal form is
%

In this paper  we consider the flow
%more general case 
of a mixture of $n$ components.
%The novelty of our contribution is the 
We study the synchronization of the velocities
in a mixture flow 
that is governed by a $2n \times 2n$ system of coupled isothermal Euler equations.
The  paper has the following structure.
%We start with the discussion of the case of gas mixtures with two components. 
%In Section \ref{isothermalEuler} we introduce
%the 
%quasi-linear 
%isothermal Euler equations
%for blended flow of two components 
%with an additional interaction term 
%as a model for hydrogen blending.
%
%In Section  \ref{sec:Riemann} we present the
%corresponding Riemann invariants
%that allow to  transform the system in diagonal form.
%We discuss a result about the well-posedness of the model.

%In Section  \ref{Synchronization2}
%we define an appropriate Lypunov function that
%allows to prove the result about exponential 
%synchronization of the two velocities.
%Our choice of the Lyapunov function is
%motivated by relative energy techniques.

In Section \ref{mixturen}, the model for 
a gas mixture that contains $n$ components is introduced.
The flow of each component is governed by 
the 
%quasi-linear 
isothermal Euler equations
%for blended flow of two components 
with an additional 
%interaction term 
velocity coupling term. 

In Section  \ref{sec:Riemannmixture} 
we show that for the  $2n \times 2n$ system, Riemann invariants
do exist that  allow to  transform the system in diagonal form.
We discuss a result about the well-posedness of the model
in the sense of semi-global classical solutions.

%as a model for hydrogen blending.

In  Section \ref{hamiltonn},
the Hamiltonian form of the system is
presented. 
Then  in Section  \ref{Lyapunovmix}, an
 appropriate 
Lyapunov function 
%for this general case 
is introduced,  that
allows to prove the result about exponential 
synchronization of the velocities.
Our choice of the Lyapunov function is
motivated by relative energy techniques.
The evolution of the Lyapunov function is analyzed in  Section \ref{evolutionmix}.
This leads to the synchronization result that is stated
in Section \ref{sec:sync} in Theorem \ref{thm7.1}.
%We present a  class of boundary conditions for which the
%velocities synchronize exponentially fast. 

%In Section \ref{sec:networks} we present
%the node conditions that model the flow 
%of the gas mixture through
%the junctions in  a gas  pipeline network.
%

%%%%%%%%%%%%%%%%%%%%%%%%%%%%%%%%%%%%%%%%%%%%%%%%%%%%

%%%%%%%%%%%%%%%%%%%%%%%%%%%%%%%%%%%%%%%%%%%%%%%%%%%%%%%%%%%%%%%%%%%%%%%%%%%%%%%%%%%%%%%%%%%%%%%%%%%%%%%%%%%%%%%%%%%%%%%%%%%%%%

%%%%%%%%%%%%%%%%%%%%%%%%%%%%%%%%%%%%%%%%%%%%%%%%%%%%%%%%%%%%%%%%%%%%%%%%%%%%%%%%%%%%%%%%%%%%%%%%%%%%%%%%%%%%%%%%%%%%%%%%%%%%%%

\section{A mixture of $n$ components: The isothermal Euler equations}
\label{mixturen}
\noindent
The one-dimensional isothermal Euler equations are a well-known 
 model for the flow of gas  through pipelines,
 see for example   \cite{bandahertyklar}.
Some  special intransient solutions 
of the  isothermal Euler equations  
are presented in \cite{iso}.
 Here we extend this model to the case of
the flow of  blended gas.

Let a natural number $n \in \{1,2,3, \ldots\}$ be given. 
%In this section 
We consider the pipeline flow of a gas 
%generalize our analysis to the case of a
mixture that consists of $n$ components.
The system that we consider is related to the study  in \cite{zlotnik}, 
but here we consider a model that consists of $n$
coupled isothermal Euler equations.

For $i \in \{1,...,n\}$,
let 
$\rho^i$ denote the
density of component $i$
and
$q^i$ the corresponding mass flow rate.
For $i\in \{1,..,n\}$, the conservation of mass is modeled by
the continuity equation 
\begin{equation}
    \label{massi}
 \rho^i_t +  q^i_x = 0.
\end{equation}
We assume  that the pressure $p^i= p_i(\rho^i)$ is given as a monotone increasing function of the density.
We  assume that this monotonicity is strict.
Examples are the isentropic gas law $p( \rho)= a \,  \rho^\gamma$ with $a >0$, $\gamma > 1$ and the model of the American Gas Association (AGA), see \cite{gugatulbrich}.
A pipe of length $L>0$ is parameterized by the interval $[0,\, L]$.
The diameter of the pipe is denoted by  $D>0$  and the  friction coefficient by  $\lambda_{fric}>0$.
Let 
%$s^e_{lope}(x) =\sin(\varphi^e(x)\,)$ and
$\theta := \frac{\lambda_{fric}}{D}$.
For $i \in \{1,...,n\}$, define the velocities 
\[v^i = \frac{q^i}{\rho^i}\]
and let
\[\rho =  \sum_{j=1}^n \rho^j ,
\;
p =  \sum_{j=1}^n p^j
. \]
Define  the barycentric velocity (see \cite{bothe2015continuum}) 
\begin{equation}
\label{vdefinition}
v = \sum_{i=1}^n  \frac{ \rho^i   }{  \rho} \, v^i. 
%\left( \frac{ \rho^i   }{  \sum_{j=1}^n \rho^j  } \right) \, v^i   
\end{equation}
For the difference
of the $i$-th velocity $v^i$ 
to the barycentric velocity we 
%introduce 
use the notation  
\begin{equation}
    \label{uidefinition}
u^i = v^i - v.
\end{equation}
%Let  real numbers $M>0$ and $\beta >0$ be given.
%With the notation $s^e= sign$
Define the coupling terms
\begin{equation}
    \label{gdefinitioni}
\bar G^i = 
\bar \Omega
%\left( \beta  + M
 %\frac{1}{2} \theta%{\rho} \,{|{v}|} 
%+ N   \right)
\,  \rho^i \,   u^i,  \; (i\in \{1,\ldots,n\})
\end{equation}
%$\omikron$
where the constant
$\bar \Omega$
%$\beta + M + N $ 
is chosen sufficiently large.
We will give precise conditions below.
%for example
%$
%N \geq  \frac{1}{2} \theta \,{|{v}|}  $.)
The flow rate of component $i$ is governed by the equation
\begin{equation}
    \label{momenti}
q^i_t+ \left(  p^i + \frac{(q^i)^2}{{  \rho}^i} \right)_x = 
- \frac{1}{2} \theta \, \frac{{q}^i \; |{q}^i|}{{\rho}^i}  -\bar  G^i´.
\end{equation}
If $\bar  G^i$ is canceled,
equations
(\ref{massi}), (\ref{momenti}) form the isothermal Euler equations,
 see e.g. \cite{bandahertyklar}.
 The interaction terms  $\bar  G^i$ model the 
 interactions between the different components of the mixture.

Similarly as in  \cite{herbertegger}, \cite{{GiesselmannGugatKunkel}} assume that 
for $i \in \{1,...,n\}$  
 smooth and 
strictly 
convex pressure potentials
\[P_i = P_i(\rho^i)\]
are given that are  connected to the pressure law by 
$ p_i'(\rho^i) = \rho^i \,  P_i''(\rho^i)$.
%$\hat p'(\rho) = \rho \, \hat{} P''(\rho)$.

\begin{example}
Similarly as in Example 1 in \cite{iso},
let a real number $\rho_0 > 0$ be given. 
For $t \geq  0$ and 
$x \in (- \infty, \infty)$ 
and $i \in \{1,  \ldots, n\}$ 
define 
$\rho^i(t, x) = \rho_0^i$.
For  a real parameter $P_0 > 0$, define
\[q^i(t, x) = \frac{1}{
\frac{P_0}{\rho_0^i} + \frac{1}{2} 
\frac{\theta}{\rho_0^i}
 t
}.
\]
Then (\ref{massi}) holds and 
we have
\[ v  \, \sum_{i=1}^n \rho_0^i = 
\sum_{i=1}^n \rho^i \, v^i
=
\sum_{i=1}^n q^i
=
\sum_{i=1}^n
\rho_0^i
\frac{1}{{P_0} + \frac{1}{2} 
{\theta}\, 
 t
}
.
\]
Thus we have
\[v = \frac{1}{P_0 + \frac{1}{2} 
{\theta}}
\]
and for all $i \in \{1, \ldots,n\}$
this yields 
$v^i = \frac{q^i}{\rho_0^i} = v$,
hence $\bar G^i =0$
and (\ref{momenti}) also holds.

\end{example}
\begin{example}
Similarly as in Example 2 in \cite{iso}, we also obtain traveling waves.
Assume that for $i \in \{1, \ldots, n\}$ we have 
$p^i(\rho^i) = a^2_i \, \rho^i$
with real numbers $a_1$, $a_2, \ldots,a_n$.
Let $C_1>0,  \ldots, C_n>0$ and $\lambda >0$ be given.
Define 
$\alpha_i(z)
=
C_i \,  \exp\left(
\frac{\lambda^2\, \theta  }
{2 \, a^2_i} \, z
\right).
$
Let
$q^i(t, \,x) = \lambda \, \alpha_i( \lambda \, t - x), \; \; 
\rho^i(t, \, x) = \alpha_i( \lambda \, t - x).
$
Then (\ref{massi}) holds and 
we have
\[ \sum_{i=1}^n \rho^i \, v^i
=
\sum_{i=1}^n q^i
= \lambda  \sum_{i=1}^n  \alpha_i( v \, t - x)
 = \lambda \sum_{i=1}^n \rho^i.
\]
This yields
$v =  \sum_{i=1}^n \frac{\rho^i}{\rho} \, v^i  =  \lambda= v^i$.
Hence $\bar G^i =0$ and (\ref{momenti}) also holds.

\end{example}
%%%%%%%%%%%%%%%%%%%%%%%%%%%%%%%%%%%%%%%%%%%%%%%%%%%%%%%%%%%%%%%%%%%%%%%%%%%%%%%%%%%%%%%%%%%%%%%%%%%%%%%%%%%%%%%%%%%%%%%%%%%%%%
\section{The system in diagonal form: Existence of semi-global classical solutions}
\label{sec:Riemannmixture} 

As stated in \cite[Chapter 7.3]{{Dafermos_book}} for   $2 \times 2$ systems of hyperbolic conservation laws
(for example for  (\ref{massi}) and (\ref{momenti}) with $\bar G_i$ replaced by zero)
we can find a system of Riemann invariants.
%Hence  we have  Riemann invariants for
%the coupled system 
%that consists of (\ref{massi}), (\ref{momenti})  for all $i\in \{1,\ldots,n\}$.
%
For the coupled system 
that consists of (\ref{massi}), (\ref{momenti})  for all $i\in \{1,\ldots,n\}$,
which  is a $2n\times 2n$ system,
we can still find Riemann invariants 
that allow to transform  the system in diagonal form.
{\mar 
The Riemann invariants can be obtained as 
potentials of suitably scaled left-hand side eigenvectors of 
the system matrix $A^Q$ of the equations in the quasilinear form
$z_t + A^Q(z) \, z_x = F^Q(z)$.
}
{The general definition of Riemann invariants can be found in 
%of systems of hyperbolic balance laws we refer to 
\cite[Sec 7.3]{Dafermos_book}.}
%For this purpose 
%We introduce 
With the notation  
$\;\tilde R(\rho)=    \int_1^\rho \frac{\sqrt{ p'(r)} }{r}   \, dr
$
we obtain the Riemann invariants
\[
 R_\pm^i (q^i, \,\rho^i)
=  \tilde R(\rho^i) \pm
\frac{q^i}{\rho^i} 
.
\]

%$ \bar R_\pm (\bar q^e, \,\bar \rho^e)
%=  \tilde R(\bar \rho^e) \pm
%\frac{\bar q^e}{\bar \rho^e} 
%$.
We have 
\begin{equation}
    \label{virepresentation}
    v^i = \frac{R^i_+ - R^i_-}{2}´, \; \rho^i =  \tilde R^{-1} \left( \frac{R^i_+ + R^i_-}{2} \right) 
\end{equation}
and
%$\rho^i =  \tilde R^{-1} \left( \frac{R^i_+ + R^i_-}{2} \right) $ 
where $  \tilde R^{-1} $ denotes the inverse function of $  \tilde R $.
We assume that as in
 the case of ideal gas, $   \tilde R^{-1} $ only
 attains values in $(0, \, \infty)$.  
This yields
\[v = \frac{1}{ \rho} \sum_{i=1}^n \rho^i \, v^i
=\frac{
\sum_{i=1}^n \left(\frac{R^i_+ - R^i_-}{2}\right)  \,  \tilde R^{-1} \left( \frac{R^i_+ + R^i_-}{2} \right) 
}{
\sum_{j=1}^n  \tilde R^{-1} \left( \frac{R^j_+ + R^j_-}{2} \right)
} .
\]
Thus we have
\[u^i = v^i - v =
 \frac{1}{ \rho} \sum_{j=1}^n \rho^j (v^i - v^j)
 =
 \frac{
\sum_{j=1}^n \left(  \frac{R^i_+ - R^i_-}{2}   - \frac{R^j_+ - R^j_-}{2} \right) \,
\tilde R^{-1} \left( \frac{R^j_+ + R^j_-}{2} \right) 
}{
\sum_{k=1}^n  \tilde R^{-1} \left( \frac{R^k_+ + R^k_-}{2} \right)
} 
.
\]

The eigenvalues of  the system matrix in
quasilinear form 
%$A(V^e)$
are
\[
 \lambda_\pm^i = 
  \frac{q^i}{\rho^i}  \pm \sqrt{ p'\rev{(\rho^i)} }, \; \; i \in \{1,\ldots, N\}.
\]

%$\bar \lambda_\pm^e = \frac{\bar q^e}{\bar \rho^e}  \pm \sqrt{ p'\rev{(\bar \rho^e)} }$
%and 
%$\hat \lambda_\pm^e = \frac{\hat q^e}{\hat \rho^e}  \pm \sqrt{ p'\rev{(\hat \rho^e)} }$.
%$\lambda_0^e = W(\hat q^e, \, \hat \rho^e)$.
{\mar
%Define $\vartheta^e(R_+^e,\, R_-^e)=  \tfrac{\theta^e}{8} |R_+^e - R_-^e| \, (R_+^e - R_-^e)  $.
Using the notation 
\begin{equation}
R^e=
\begin{pmatrix}
R_+^1
\\
 R_-^1
\\
 R_+^2
\\
 R_-^2
\\
\vdots
\\
 R_+^N
\\
 R_-^N
     \end{pmatrix},
     \;
     D(R^e) = 
\begin{pmatrix}
\lambda_+^1
%(R^e_+, \, R^e_-) 
&  0 & 0 & 0 & \ldots & 0
\\
%P_{  \rho_{(h)}}
0 &  \lambda_-^1 & 0 
%(R^e_+, \, R^e_-)  
& 0 &   \ldots & 0
\\ 
0  &   0  &  \lambda^2_+
& 0 &   \ldots & 0
\\
0 & 0 &  0 & \ddots &  &     0
\\ 
0 & 0 &  \ldots & 0 & \lambda_+^N &  0
\\
0 & 0 &    \ldots & 0  & 0 & \lambda_-^N 
%(R^e_+, \, R^e_-)
\end{pmatrix}
,
\,
%S^e(R^e) = 
%\left(
%\begin{array}{c}
%\begin{pmatrix}
%- \vartheta^e(R_+^e,\, R_-^e)
%\tfrac{\theta^e}{8} |R_+^e - R_-^e| \, (R_+^e - R_-^e)
%\\
%\vartheta^e(R_+^e,\, R_-^e)
% \tfrac{\theta^e}{8}   |R_+^e - R_-^e| \, (R_+^e - R_-^e)
%\\
%0
%\end{pmatrix}
%\end{array}
%\right)
\end{equation}
(with $R^e \in {\mathbb R}^{2n}$)
we can write system  (\ref{massi}), (\ref{momenti}), ($i\in \{1,\ldots,n\}$),
in the diagonal form 
%\begin{multline}
\begin{equation}
   \label{hydrodiag}
%
%system (\ref{isothermaleuler}), (\ref{hydrogen})
%can be written in the quasi-linear form
R^e_t + D(R^e) \, R^e_x = S^e(R^e)
\end{equation}
%\end{multline}
where $S^e$ denotes the corresponding source term.

}
%
%
%
%we can write the system
%\eqref{isothermaleuler}, \eqref{hydrogen}
%in the diagonal form 
%\begin{multline}
%  %  \label{hydrodiag}
%%\left\{
%\begin{pmatrix}
%R_+^e
%\\
%R_-^e
%\\
%R_0^e 
%%c_{(h)} 
%\end{pmatrix}_t
%+
%\begin{pmatrix}
%\lambda_+(R_+^e, R_-^e) & 0 & 0
%\\
%%\partial_{\hat \rho} 
% 0 & \lambda_-(R_+^e, R_-^e) &
%%\partial_{  \rho_{(h)}}
%0
%\\
% 0 &   0 & \lambda_0(R_+^e, \, R_-^e)
%\end{pmatrix}
%\begin{pmatrix}
%R_+^e
%\\
%R_-^e
%\\
%R_0^e 
%%c_{(h)} 
%\end{pmatrix}_x
%\\
%=
%\begin{pmatrix}
%- \tfrac{\theta^e}{8} |R_+^e - R_-^e| \, (R_+^e - R_-^e)
%\\
% \tfrac{\theta^e}{8}   |R_+^e - R_-^e| \, (R_+^e - R_-^e)
%\\
%0
%\end{pmatrix}
%.
%\end{multline}
Here the eigenvalues are represented as functions of
the Riemann invariants, i.e. 
%for example
%\begin{eqnarray}
%\label{lambdaminusgleichung}
%\tilde \lambda^e_-  &  = &
%\sqrt{R^e_s\, T^e}\,
%\left[
%\frac{R_+^e - R_-^e}{2} - 1 - \alpha^e\, \exp\left(\frac{R_+^e + R_-^e}{2}\right)
%\right],
%\\
%\label{lambdaplusgleichung}
%\tilde \lambda^e_\pm
%&
%=
%&
%\sqrt{R^e_s\, T^e}\,
%\left[
%\frac{R_+^e - R_-^e}{2} + 1 + \alpha^e\, \exp\left(\frac{R_+^e + R_-^e}{2}\right)
%\right].
\[
%\label{lambdapmgleichung}
 \lambda^i_\pm
%&
=
%&
\frac{R_+^i - R_-^i}{2} 
\pm 
\sqrt{p'\left(\tilde R^{-1}   \left(\frac{R_+^i + R_-^i}{2}\right) \right)  }
.
\]
%\end{eqnarray}
%Specific formulas for Riemann invariants  for the isentropic law $p(\rho)=a %\rho^\gamma$ can be found in \cite[Chapter 7.3]{dafermosbook}
%while 
The case for the AGA model for real gas  is presented in detail in  \cite{gugatulbrich}.
%Note that the Riemann invariants are not uniquely determined.
 %up to operations that leave the direction of the gradient unchanged.
%for $e\in E$ the Riemann invariants  $R^e_-$, $R^e_+$  of the system are given by
%\[R_-^e(p^e,\, q^e)
%=
%{\rm ln }(p^e)  - \sqrt{R^e\, T^e}\, \frac{q^e}{p^e} - \alpha^e \, \sqrt{R^e\, T^e}\, q^e,
%=
%{\rm ln }(p^e)  - \sqrt{R^e_s\, T^e}\, \frac{q^e}{p^e} \left( 1  + \alpha^e \, p^e\right),
 %\sqrt{R^e\, T^e}\, q^e,
%\]
% \[R_+^e(p^e,\, q^e) =
 % {\rm ln }(p^e)  + \sqrt{R^e\, T^e}\, \frac{q^e}{p^e} + \alpha^e\, \sqrt{R^e\, T^e}\, q^e.
% {\rm ln }(p^e)  + \sqrt{R^e_s\, T^e}\, \frac{q^e}{p^e} \left( 1  + \alpha^e \, p^e\right).
% \]
%%%%%%%%%%%%%%%%%%%%%%%%%%%%%%%%%%%%%%%%
For $ i \in \{1, \ldots, n\}$  define 
\[\vartheta_i(R^e)=  \tfrac{\theta}{8} |R_+^i - R_-^i| \, (R_+^i - R_-^i)
+
 \bar \Omega \, \tilde R^{-1} \left( \frac{R^i_+ + R^i_-}{2} \right) \, u^i(R^e)
. \]
Then we have
\[
S^e(R^e) = 
\begin{pmatrix}
%\left(
%\begin{array}{c}
- \vartheta_1(R^e)
%\tfrac{\theta^e}{8} |R_+^e - R_-^e| \, (R_+^e - R_-^e)
\\
\vartheta_1(R^e)
% \tfrac{\theta^e}{8}   |R_+^e - R_-^e| \, (R_+^e - R_-^e)
\\
\vdots
\\
- \vartheta_n(R^e)
%\tfrac{\theta^e}{8} |R_+^e - R_-^e| \, (R_+^e - R_-^e)
\\
\vartheta_n(R^e)
%\end{array}
%\right)
\end{pmatrix}
.
\]
%{
%\mar
%Note that 
%the transformation back to the physical variables yields density values $\hat \rho^e >0$.
%in the diagonal form in  the source term
%in contrast to  (\ref{isothermaleuler})
%division 
%a quotient does not occur.
%by $\hat \rho^e$ 
%{isothermaleuler}
%}
The diagonal form of the system is of
interest both from the theoretical and the numerical point of view:
Similarly as in \cite{gugatulbrich}, it  allows to show the existence of a solution in 
the sense of the characteristics since
the system can be rewritten as a system of integral equations along the
characteristic curves.
This solution can be characterized as the unique solution
of an appropriately defined fixed-point iteration.
On the other hand, assuming that the signs of  the eigenvalues 
are constant,
%in this situation 
for the numerical approximation of the system state
an upwind-downwind discretization  makes sense,
see \cite{grundel}.
This  discretization  mimics the flow of information along
the characteristic curves whose
direction is determined by the signs of the eigenvalues.

%%%%%%%%%%%%%%%%%%%%%%%%%%%%%%%%%%%%%%%%%%%%%%%%%%%%%%%%%%%%%%%%%%%%%%%%%%%%%%%%%%%%%%%%%%%%

%%%%%%%%%%%%%%%%%%%%%%%%%%%%%%%%%%%%%%%%%%%%%%%%%%%%%%%%%%%%%%%%%%%%%%%%%%%%%%%%%%%%%%%%%%%%%%%%%%%%%%%%%%%%%%%%%%%%%
\subsection{Well-posedness}
To obtain a well-posed problem, we 
complete the system with initial- and boundary data.
For $i\in \{1,\ldots,n\}$ let  initial velocities
$v^i_0 \in C^1([0, \, L])$
and initial densities
$\rho^i_0  \in C^1([0, \, L])$
be given.

With the initial conditions
\begin{equation}
\label{initial}
v^i(0, \, x) = v^i_0,
\;
\rho^i(0, \, x ) = \rho^i_0, \; x\in [0, \, L], \; i\in \{1,\ldots,n\}
\end{equation}
also the values of the Riemann invariants at the time $t=0$
are prescribed.

In addition to the initial condition
(\ref{initial}) also boundary conditions 
at $x=0$ and $x=L$ have to be prescribed.
At $x_B=0$ assume that 
\begin{equation}
\label{5.10}
v^i(t, \, x_B) = \bar v(t) ,\;   i \in \{1,\ldots,n\}
\end{equation}
were $\bar v(t) > 0 $ is a continuously differentiable function.
Due to (\ref{virepresentation}), we can state
(\ref{5.10}) in terms of Riemann invariants as
\begin{equation}
\label{velocityriemann}
R^i_+(t, \, x_B) = R^i_-(t, \, x_B) + 2 \, \bar v(t) ,\;   i \in \{1,\ldots,n\}.
\end{equation}

At $x = L,$
for $i \in \{1, \ldots,n\}$ 
we consider the boundary conditions
%that are given in terms  of Riemann-invariants as 
%\begin{equation}
%\label{riemannx=0}
%R^i_-(t, \, L) = \Psi^i(t, \,   R^1_+(t, \, L),\ldots, \, R^n_+(t, %\, L))
%\end{equation}
%with continuously differentiable maps $\Psi^i$. 
%
%As an example in terms of the physical variables,
% apart from the adaptation of (\ref{5.10}) to $x_B=0$ 
%consider the boundary conditions
%of  feedback law
\begin{equation}
\label{feedbackrho}
\rho^i(t, \, L)   = \bar \rho^i(t) ,\;   i \in \{1,\ldots,n\}
%\Phi^i (  t, \, v^1(t, \, x_B) ,  \, \ldots, v^n(t, \, x_B)) % , \, I(t)) 
\end{equation}
where the maps $\bar \rho^i$ are continuously differentiable and attain values in $[0, \infty)$.
In terms of the Riemann invariants,
due to
(\ref{virepresentation}) 
this is equivalent to
\begin{equation}
\label{densityriemann}
R^i_-(t, \, L) = - R^i_+(t, \, L) + 2 \,  \tilde R( \bar \rho(t)).
\end{equation}

%or
%\begin{equation}
%\label{feedback}
%v^i(t, \, x_B )   = \Phi^i (  t, \, \rho^1(t, \, x_B) ,  \, \ldots, \rho^n(t, \, x_B)) % , \, I(t)) 
%\end{equation}

%where the maps $\Phi^i$ are continuously differentiable
%and attain values in $[0, \infty)$
%in  the case of (\ref{feedback}), 
%and strictly positive values  for  (\ref{feedbackrho}). 

%Note that if for all $i, j \in \{1,\ldots,n\}$ we have
%$ p_i'(\rho^i) = \rho^i \,  P_i''(\rho^i) =  p_j'(\rho^j)$
%(which is the case if the gas laws for the components  are identical)
%we have $I(t)\not=0$ in general 

We assume that the $C^1$-compatibility conditions holds.
For (\ref{5.10}) and (\ref{feedbackrho}) this means 
that $\bar v(0)$, $\bar v_t(0)$ and
$\bar \rho^i(0)$, $\bar \rho^i_t(0)$ 
are  compatible with the values of 
$v^i(0, x_B)$, $v^i_t(0, \, x_B)$ and  $\rho^i(0, \, L)$,  $\rho^i_t(0, \, L) $
% and $I'(0)$
that  are obtained from the initial data  and  the partial differential equation.
%that is the $C^1$-compatibility conditions holds.
%For (\ref{5.10}) and (\ref{feedbackrho}) this means 
%for example
To be precise   we have the conditions 
\[v^i(0, \, x_B)    = \bar v(0), \, \rho^i(0, \, L) = \bar \rho^i(0), 
%\Phi^i (  x_B, \, \rho^1_0 ,  \, \ldots, \rho^n_0)
%, I(0) ) 
\]
\[
\bar \rho^i_t(0) =  - \partial_x( v^i(0, \, \cdot) \, \rho^i(0, \cdot)   )|_{x=L}
\]
and
\[
  \bar v_t(0) 
  = 
  -\partial_x 
  \left.
  \left( \frac{1}{2} v^2(0,\, \cdot) + P_i'(\rho^i( 0, \cdot)) \right)\right|_{x=x_B} 
   - 
   \left[\frac{1}{2}\theta \,   v^i
   \, |v^i|
+
%(\beta + M + M ) 
\bar\Omega 
\,    ( v^i - v )
\right](0, x_B)
%   
%(\bar \rho^i \, \bar v^i)_t
%+ \left(  p^i(0, \cdot) + \frac{(q^i(0, \cdot))^2}{{  \rho}^i(0, \cdot)} \right)_x|_{x=0}  = 
%- \frac{1}{2} \theta \, \frac{{q}^i \; |{q}^i|}{{\rho}^i}  -\bar  G^i´.
.
\]
For the derivation of the last condition, see also
Section \ref{hamiltonn}  on the system in Hamiltonian form below.

%in the case (\ref{feedback})
%and that for $t=0$ also  the equations 
%\[v^i_t(t, \, x_B)   = \partial_t \Phi^i (  t, \, \rho^1(t, \, x_B) ,  \, \ldots, \rho^n(t, \, x_B)   ) %, I(t)  )
%\]
%\[
%+ \sum_{i=1}^n \partial_{\rho^i} \Phi^i (  t, \, \rho^1(t, \, x_B) ,  \, \ldots, \rho^n(t, \, x_B) ) % , I(t) )
%\rho^i_t(t,x_B)
%\]
%\[
%+
%\partial_{I }\Phi^i (  t, \, \rho^1(t, \, 0) ,  \, \ldots, \rho^n(t, \, 0) , I(t) ) I'(t)
%\]
%
%are  compatible with the values of 
%$v^i_t(x_B, \, 0 )$ and  $ \rho^i_t(x_B ,0) $ % and $I'(0)$
%that  are obtained from the initial data  and  the partial differential equation,
%that is the $C^1$-compatibility conditions holds.

Then the  well-posedness of system 
(\ref{initial}), 
(\ref{5.10}), 
(\ref{feedbackrho}), 
(\ref{massi}), (\ref{momenti}) 
%
%(\ref{hydrodiag})
follows from the theory of semi-global classical solutions for quasilinear hyperbolic systems,
see for example 
\cite{li2001semi},
\cite{wang2006exact} 
or the  results in \cite{{GiesselmannGugatKunkel}}.
The theory of semi-global classical solutions 
asserts the following:
Given a classical intransient solution $\bar z$ of the system
that is $C^1$-compatible with the boundary conditions
%feedback law (\ref{feedback}) 
and a time horizon $T$, there exists a number $\varepsilon(T)>0$ 
with the following property:
For initial data in a $C^1$-neigbourhood of $\bar z$ with radius 
$\varepsilon(T)$ 
that are  $C^1$-compatible with the boundary conditions  
%(\ref{feedback}) 
%and $C^1$-compatible  boundary data in
%a $C^1$-neigbourhood with radius  $\varepsilon(T)$  of the constant boundary data
%corresponding to $\bar z$ on
%the time interval $[0, T]$
there exists a unique classical solution
on the space-time rectangle $[0, \, T] \times [0, \, L]$. 
This solution  satisfies an a-priori inequality
that yields an upper bound for the
$C^1$-norm of the solution on $[0, \, T] \times [0, \, L]$
that is linear in terms of
the $C^1$-norms of the initial data
%and the boundary data.
(see also \cite{bastincoron}, for the existence of solutions with
$H^2$-regularity).
This  a-priori inequality implies in particular that
by choosing 
$\varepsilon(T)>0$ 
sufficiently small, 
it can be guaranteed that 
the solutions satisfies  
prescribed $C^1$-bounds.

Note that system 
(\ref{initial}), 
(\ref{5.10}), 
(\ref{feedbackrho}), 
(\ref{massi}), (\ref{momenti}) 
with $\bar v(t)=0$ 
and $\bar \rho^i(t)  = \tilde R^{-1}(J^i) $
has constant solutions
$R^i_\pm = J^i$ for all $i \in \{1, \ldots, n\}$.
We study the classical solutions in a $C^1$-neighbourhood of a constant solution.
%%%%%%%%%%%%%%%%%%%%%%%%%%%%%%%%%%%%%%%%%%%%%%%%%%%%%%%%%%%%%%%%%%%%%%%%%%%%%%%%%%%%%%%%%%%%%%%%%%%%%%%%%%%%%%%%%%%

First we state the results for the system in diagonal form.
In this problem we have 
the initial conditions for the Riemann invariants
\begin{equation}
\label{initialriemann}
R^i_\pm (0, \, x)  = R^{(i, 0)}_\pm(x)
\end{equation}
with continuously differentiable functions 
$  R^{(i, 0)}_\pm $.

Then we have the initial boundary -value problem problem
(\ref{initialriemann}), 
(\ref{velocityriemann}), 
(\ref{densityriemann}),
(\ref{hydrodiag}).

\begin{theorem}[see for example \cite{wang2006exact}]
\label{semiglobal}
% [Semi-global classical solutions ({\sc Ta-Tsien Li} (Fudan university, Shanghai),...)]
Let $T>0$
and  for $i \in \{1, \, \ldots,n\}$  let  real numbers $J^i$
%and a 
%sufficiently small 
%number $M>0$ 
be given.
Then there exists 
%is a number 
$\varepsilon(T)
%, \, M)
>0$ such that
for all {initial data}
$  R^{(i, 0)}_\pm   \in C^{1}(0,\, L)$
%($e\in E$)
such that
\[
\tilde{I} :=
\max_{ i \in \{1, \, \ldots,n\}   }
\|
 R^{(i, 0)}_\pm  
-J^i\|_{ C^1(0,\, L) } \leq \varepsilon(T)
%,\, M)
\]
and all functions
$\bar v(t)$, $\bar \rho^i(t)\in C^1(0,\, T)$ 
such that
\[
\tilde {C} :=
\max_{ i \in \{1, \, \ldots,n\}   }
\{\|\bar v\|_{C^1(0, T)}, \;    \|\bar \rho^i - \tilde R^{-1}( J^i)   \|_{C^1(0, T)} 
\}
\leq \varepsilon(T)
%,\, M) 
\]
that satisfy the  \emph{$C^1$-compatibility conditions} for 
%${\bf (S) }$,  
(\ref{hydrodiag}), 
(\ref{velocityriemann}), 
(\ref{densityriemann}),
(\ref{initialriemann})
there exists a
unique classical solution  of
(\ref{hydrodiag}), 
(\ref{velocityriemann}), 
(\ref{densityriemann}),
(\ref{initialriemann})
 on $[0, \, T]$.

% that satisfies the integral equations
% (\ref{integralequation}) for all $e\in E$
% along the characteristic curves 
% with
% %$S^e_{+}$,  $S^e_{-} \in L^\infty((0,\, L^e)\times (0,\, T))$
% $S^e_{+}$,  $S^e_{-} \in C^1((0,\, T)\times (0,\, L^e))$
% ($e\in E$)
% and the boundary condition
% (\ref{rbriemann2neu}) 
% at the boundary nodes 
% and the node condition
% (\ref{Omegav}) at the interior nodes 
 % almost everywhere 
% in $[0, \, T]$ 
%
%such that for all $e\in E$ we have
%\begin{equation}
%\label{lunendlichbound}
%\|S^e_{\pm} - J\|_{C^1( (0,\, T) \times (0,\, L^e) ) } \leq M.
%\end{equation}
%
%Along the characteristic curves, the solution
%is absolutely continuous.
%This solution depends in a stable way  on the initial and boundary data in the sense that
%for initial data
%$ 
%\|z^e_{\pm}-J\|_{C^1(0, L^e)} \leq \varepsilon(T,\, M)
%$
%and control functions
%$v^e \in C^1(0,\, T)$  ($e\in E$)
%that satisfy the  $C^1$ compatibility conditions for 
%${\bf (S) }$  
%such that
%$
%\|v^e-J \|_{L^\infty(0, T)} \leq \varepsilon(T,\, M)
%$
%for the corresponding solution $\tilde S^e_\pm$ we have the inequality
%\[
%\| \tilde S^e_\pm - S^e_{\pm}\|_{C^1( (0,\, T) \times (0,\, L^e) ) } \leq 
%C(T)\; \max_{e \in E} \{ \|y^e_{\pm} - z^e_{\pm}\|_{C^1(0, L^e)},\,  \|u^e - v^e \|_{C^1(0, T)}\}
%\]
%where $C(T)> 0$ is a constant that does not depend on $z^e_{\pm}$ or $v^e$.
%

There exists a constant $\hat C_T>0$ such that 
the solution satisfies the 
\textbf{a  priori bound}
\begin{equation}
    \label{aprioribound2020}
 %\max_{s \in [0,\, T]}
 \max_{ i \in \{1, \, \ldots,n\}   }
% \left\{
% \|R^e_+(s,\, x) - J\|_{L^\infty(0, \, L^e)},
% \;
% \|R^e_-(s,\,x)- J \|_{L^\infty(0, \, L^e)}
% \right\}
\left\{
\|R^i_\pm - J^i\|_{C^1( (0,\, T) \times (0,\, L) ) }
%{C^1(0, \, L^e)}
\right\}
%\[
\leq
\hat C_T \,
\max\{\tilde I, \, \tilde C\}
%\,  \varepsilon(T,\, M) 
%3\, \varepsilon(T,\, M) \exp( 48  \, \max_{e\in E} \nu^e\, M \,T)
.
\end{equation}

\end{theorem}
%%%%%%%%%%%%%%%%%%%%%%%%%%%%%%%%%%%%%%%%%%%%%%%%%%%%%%%%%%%%%%%%%%%%%%%%%%%%%%%%%%%%%%%%%%%%%%%%%%%%%%%%%%%%%%%%%%%%%%%%%%

Theorem  \ref{semiglobal} yields the following result in terms of the physical variables.

\begin{theorem}%[see for example \cite{wang2006exact}]
\label{semiglobalphysical}
% [Semi-global classical solutions ({\sc Ta-Tsien Li} (Fudan university, Shanghai),...)]
Let $T>0$  and  for $i \in \{1, \, \ldots,n\}$  let  real numbers $J^i$
%and a 
%sufficiently small 
%number $M>0$ 
be given.
Then there exists 
%is a number 
$\varepsilon(T)
%, \, M)
>0$ such that
for all {initial data}
$  v^i_0$, $\rho^i_0   \in C^{1}(0,\, L)$
%($e\in E$)
such that
\[
\tilde{I} :=
\max_{ i \in \{1, \, \ldots,n\}   }
\|
v^i_0
\|_{ C^1(0,\, L) } \leq \varepsilon(T),
\;
 \|\rho^i_0  - \tilde R^{-1}( J^i)   \|_{C^1(0, T)} \leq \varepsilon(T),
%,\, M)
\]
and all functions
$\bar v(t)$, $\bar \rho^i(t)\in C^1(0,\, T)$ 
such that
\[
\tilde {C} :=
\max_{ i \in \{1, \, \ldots,n\}   }
\{\|\bar v\|_{C^1(0, T)}, \;    \|\bar \rho^i - \tilde R^{-1}( J^i)   \|_{C^1(0, T)} 
\}
\leq \varepsilon(T)
%,\, M) 
\]
that satisfy the  \emph{$C^1$-compatibility conditions} for 
%${\bf (S) }$,  
(\ref{massi}), 
(\ref{momenti}), 
(\ref{initial}), 
(\ref{5.10}), 
(\ref{feedbackrho})
there exists a
unique classical solution  of
(\ref{massi}), 
(\ref{momenti}), 
(\ref{initial}), 
(\ref{5.10}), 
(\ref{feedbackrho}) 
 on $[0, \, T]$.

There exists a constant $\hat C_T>0$ such that 
the solution satisfies the 
a  priori bound 
\begin{equation}
    \label{aprioribound2025}
 %\max_{s \in [0,\, T]}
 \max_{ i \in \{1, \, \ldots,n\}   }
% \left\{
% \|R^e_+(s,\, x) - J\|_{L^\infty(0, \, L^e)},
% \;
% \|R^e_-(s,\,x)- J \|_{L^\infty(0, \, L^e)}
% \right\}
\left\{
\| v^i \|_{C^1( (0,\, T) \times (0,\, L) ) }, \;  
\|\rho^i - \tilde R^{-1} (J^i) \|_{C^1( (0,\, T) \times (0,\, L) ) }
%{C^1(0, \, L^e)}
\right\}
%\[
\leq
\hat C_T \,
\max\{\tilde I, \, \tilde C\}
%\,  \varepsilon(T,\, M) 
%3\, \varepsilon(T,\, M) \exp( 48  \, \max_{e\in E} \nu^e\, M \,T)
.
\end{equation}
%Moreover, if $\varepsilon(T)$ is chosen 
%sufficiently small there exists a constant $M_0 > 0$  such that for
%all $i \in \{1, \, \ldots,n\}  $
%and all $(t, \, x) \in  [0, \, T] \times [0, \, L]$ we have
%\begin{equation}
%\label{rhoischranke}
%    \rho^i(t, \, x)   
    %\frac{\rho^i}{2} 
%    \geq  2\,  M_0^2. 
%\end{equation}

\end{theorem}
%%%%%%%%%%%%%%%%%%%%%%%%%%%%%%%%%%%%%%%%%%%%%%%%%%%%%%%%%%%%%%%%%%%%%%%%%%%%%%%%%%%%%%%%%%%%%%%%%%%%%%%%%%%%%%%%%%%%%%%%%%

\begin{remark}
{%\color{blue}
It is well-known that for quasilinear hyperbolic systems,
classical  solutions can break down after finite time, 
since for example rarefaction fans,  shock waves,
or contact discontinuities 
can evolve even from arbitrarily smooth initial states,
see eg. \cite{rinaldo2}.
It is however important to note that
in the operation of gas pipelines,
an important aim is to control the system in 
such a way that shocks are avoided, 
since they can lead to damage of
the pipeline.
Hydrogen embrittlement is a degradation effect
that has to be taken into account as soon
as hydrogen is transported in pipelines made of steel.
The fatigue crack growth behavior is influenced by the loading conditions,
see e.g. \cite{laureys}.

From the point of view of the system operation it is of interest 
whether it is possible to avoid the evolution of singularities.
From the mathematical point of view, this
corresponds to the investigation of classical solution of
the system, and in this contribution we will restrict our attention
to this framework.
However, it is clear that for a complete analysis it is also important to
understand whether  also for  solutions with  less regularity, the
velocities still synchronize.
We will come back to this point in the conclusions.

%the influence of  source term

}

\end{remark}

%%%%%%%%%%%%%%%%%%%%%%%%%%%%%%%%%%%%%%%%%%%%%%%%%%%%%%%%%%%%%%%%%%%%%%%%%%%%%%%%%%%%%%%%%%%%%%%%%%%%%%%%%%%%%%%%%%%%%%%%%%

%%%%%%%%%%%%%%%%%%%%%%%%%%%%%%%%%%%%%%%%%%%%%%%%%%%%%%%%%%%%%%%%%%%%%%%%%%%%%%%%%%%%%%%%%%%%%%%%%%%%%%%%%%%%%%%%%%%%%%%%%%%%%%

\section{The system in Hamiltonian form}
\label{hamiltonn}
For $i \in \{1,...,n\}$  define the energy densities %functionals 
%\[h(a,b) = \frac{1}{2} a\, b^2 + p( a)\]
%and
\[H^i(a, \, b) = \frac{1}{2} a\, b^2 + P_i( a),\]
and let 
\[
H^{mix}(a_1, a_2, \ldots, a_n, \, b) = \frac{1}{2} \sum_{i=1}^n  a_i \, b^2 + \sum_{i=1}^n P_i( a_i )
=
\sum_{i=1}^n H^i(a_i, \, b) .
\]

We use the notation $
%(H^i)'=
(H^i_\rho, \, H^i_v)$ for the 
%variational 
partial derivatives with
\begin{equation}
\label{varderi}
H^i_{\rho^i} (\rho^i, \, v^i)  = \frac{1}{2} ( v^i)^2   + P_i'(\rho^i), \; 
H^i_{v^i} (\rho^i, \, v^i) = \rho^i \, v^i, 
\end{equation}
cf. (41) in  \cite{GiesselmannGugatKunkel}. 
We have 
\begin{equation}
\label{vableitung}
H^{mix}_{v}(\rho_1, \rho_2, \ldots, \rho_n, \, v) = \sum_{i=1}^n  H^i_v(\rho_i, \, v)  = \rho \, v,
\end{equation}
and 
%with a slight abuse of notation we define  
\[
H^{mix}_{\rho^i}(\rho_1, \rho_2, \ldots, \rho_n, \, v) = 
\frac{1}{2} v^2 + 
%\sum_{i=1}^n 
P_i'(\rho^i)
=
H^i_{\rho^i} (\rho^i, \, v).
\]

 The Hamiltonian form  of  system (\ref{massi}),  (\ref{momenti}),  is
\begin{equation}
 \label{hamiltonianmix}
\left\{
\begin{array}{l}
\rho^i_t + \partial_x  H^i_{v^i} (\rho^i, \, v^i)  =   0,
\\
  v^i_t + \partial_x  H^i_{\rho^i}  (\rho^i, \, v^i)    =   - \frac{1}{2}\theta \,   v^i\, |v^i|
-
%(\beta + M + M ) 
\bar\Omega 
\,    ( v^i - v)
% - 
 %\frac{1}{\rho}
% \frac{1}{2} \theta \, %s^e \,
%{\bar\rho}
%{\rho} 
%\,{|{v}|} \, (v^i - v)
.
 \end{array}
 \right.
\end{equation}

This implies in turn 
\begin{theorem}
Assume that for
 $i \in \{1,...,n\}$,
a classical  solution of
(\ref{hamiltonianmix}) is given on
$[0, T]\times [0, L]$
with $\rho >0$ on
$[0, T]\times [0, L]$.
Then for the barycentric velocity $v$ defined in
(\ref{vdefinition}) we have the 
partial differential equations 
 \begin{equation}
 \label{hamiltonianmixturn}
\left\{
\begin{array}{l}
 \rho_t + \partial_x H^{mix}_v = 0,
 \\
% v_t + \partial_x H_\rho 
% +
%  \left(  \frac{\bar\rho}{\rho}\,   \partial_x ( \bar v^e \, \bar u^e)
%+   \frac{\hat\rho}{\rho} \,   \partial_x (\hat v^e \, \hat u^e) \right) 
%+
%\frac{\bar \rho_x}{\rho} \, \bar v \, \bar u 
%+
%\frac{\hat \rho_x}{\rho} \, \hat v \, \hat u 
%\\
%=
v_t + 
v \, v_x +
\left(
\sum_{i=1}^n  
 \frac{ \rho^i }{\rho}  \,  P_i''(\rho^i)\, \rho^i_x 
\right)
=
\frac{1}{\rho} \, 
\left( 
-
\sum_{i=1}^n  
\left(
\rho^i\, (u^i)^2
\right)_x 
-
\sum_{i=1}^n  
 \frac{1}{2}\theta \, \rho^i \,   v^i\, |v^i|
 \right).
%\partial_x H^{mix}_\rho 
 \end{array}
 \right.
\end{equation}
\end{theorem}
\begin{remark}
Note that the equation
%second  equation in 
%(\ref{hamiltonianmixturn}) and
\begin{equation}
    \label{vmischung}
v_t + 
v \, v_x +
\left(
\sum_{i=1}^n  
 \frac{ \rho^i }{\rho}  \,  P_i''(\rho^i)\, \rho^i_x 
\right)
=
-
 \frac{1}{2}\theta  \,   v\, |v|
%\partial_x H^{mix}_\rho 
\end{equation}
and  the $n$ continuity equations 
$\rho^i_t + (\rho^i \, v)_x=0$
($i \in \{1,\ldots, n\}$)  
determine a state
$(\rho^1, \ldots, \rho^n, v)$
from an initial state with 
the same velocities  (e.g. $v^i = v =0$ for all $i \in \{1,\ldots, n\}$)
and compatible boundary data where also $v^i  = v$
for all $i \in \{1,\ldots, n\}$.
Equation (\ref{vmischung}) 
follows from the  second  equation in (\ref{hamiltonianmixturn}) 
in the case that we have $v^i  = v$
for all $i \in \{1,\ldots, n\}$.

System  (\ref{hamiltonianmixturn})  can be written in the form
\begin{equation}
 \label{hamiltonianmixturnconvex}
\left\{
\begin{array}{l}
 \rho_t + \sum_{i=1}^n  
 \frac{ \rho^i }{\rho}  \partial_x H^{mix}_v(\rho_1,  \ldots, \rho_n, \, v)  = 0,
 \\
v_t + 
\sum_{i=1}^n  
 \frac{ \rho^i }{\rho} 
 \partial_x H^{mix}_{\rho^i}(\rho_1,  \ldots, \rho_n, \, v) 
 =
 \frac{1}{\rho} \, 
\left( 
-
\sum_{i=1}^n  
\left(
\rho^i\, (u^i)^2
\right)_x 
-
\sum_{i=1}^n  
 \frac{1}{2}\theta \, \rho^i \,   v^i\, |v^i|
 \right)
 \end{array}
 \right.
\end{equation}
which contains a convex combination of the contributions  of the components
of the gas mixture.

\end{remark}

\begin{proof}
The first equation  in  (\ref{hamiltonianmixturn}) 
follows by summing up the first equations in (\ref{hamiltonianmix})
which are  the  continuity equations for the  $\rho^i$ and  with    (\ref{vableitung}).
%The derivation of the second equation   in  (\ref{hamiltonianmixturn})  is more involved. 
%given below.

Now we  show how the second equation in  (\ref{hamiltonianmixturn}) can be derived from  (\ref{hamiltonianmix}). 
Note that by (\ref{vdefinition}) we have 
\begin{equation}
\label{rhovgleichung}
 \rho \, v = \sum_{i=1}^n  \rho^i \, v^i.
 \end{equation}
%This implies
%\[
%\rho \, v_t = v^2 \, \rho_x + \rho \, v \, v_x  - \sum_{i=1}^n  \left( (v^i)^2 \, \rho^i_x + \rho^i \, v^i \, v^i_x    \right).
%\]
%This can be seen as follows.
%
For  $q = \rho \, v = \sum_{i=1}^n \rho^i \, v^i  = \sum_{i=1}^n q^i $
we have
\begin{equation}
    \label{vtgleichung}
    v_t = \frac{1}{\rho} \, q_t  + \frac{v}{\rho}\, q_x
=
\frac{1}{\rho} \, q_t  + v \, v_x + v^2 \, \frac{\rho_x}{\rho}.
\end{equation}
Similarly, for all $i \in \{1,\ldots, n\}$ we have 
\[v^i_t = \frac{1}{\rho^i} \, q_t^i  + \frac{v^i}{\rho^i}\, q_x^i
=
\frac{1}{\rho^i} \, q_t^i  + v^i \, v_x^i + (v^i)^2 \, \frac{\rho_x^i}{\rho^i}.
\]

Note that  due to (\ref{rhovgleichung}) we have 
\begin{equation}
    \label{rhoiuigleichung}
\sum_{i=1}^n  \rho^i \, u^i
=
\sum_{i=1}^n
 \rho^i \, (v^i - v )  = \left(\sum_{i=1}^n  \rho^i \,  v^i \right)  -v\, \sum_{i=1}^n  \rho^i 
= \rho \, v - v \rho =0.
\end{equation}
Thus we have 
\begin{equation}
    \label{uiquadratgleichung}
\sum_{i=1}^n  \rho_i \, (v^i)^2
=
\sum_{i=1}^n  \rho_i \, (v + u^i)^2
=
\sum_{i=1}^n
\left[\rho^i v^2 +  2\, v\,  \rho^i \, u^i  + \rho^i\, (u^i)^2 \right]
\end{equation}
\[
=
\rho \, v^2 + \sum_{i=1}^n \rho^i\, (u^i)^2
.
\]
Due to (\ref{gdefinitioni}) and (\ref{rhoiuigleichung})    we have 
$\sum_{i=1}^n \bar G^i =0$. Hence  we obtain  
%$q_t =  \sum_{i=1}^n \rho^i \, v^i_t + v^i \,  \rho^i_t$
%
\[
q_t 
=\sum_{i=1}^n q^i_t 
= \sum_{i=1}^n
\left[
-
\left(  p^i + \frac{(q^i)^2}{{  \rho}^i} \right)_x 
- \frac{1}{2} \theta \, \frac{{q}^i \; |{q}^i|}{{\rho}^i}  -\bar  G^i
\right] 
\]
\[
= \sum_{i=1}^n
\left[
-
  p_i'(\rho^i) \, \rho^i_x 
  -
  \left( \rho^i \, (v^i)^2\right)_x 
- \frac{1}{2} \theta \, \frac{{q}^i \; |{q}^i|}{{\rho}^i}
\right] 
\]
\[
= 
\left( \sum_{i=1}^n
- \rho^i \,  P_i''(\rho^i)\, \rho^i_x \right)
 -  \left(  \rho \, v^2 + \sum_{i=1}^n \rho^i\, (u^i)^2 \right)_x
-  \sum_{i=1}^n \frac{1}{2} \theta \,  {{\rho}^i} \,  {{v}^i \; |{v}^i|}.
\]
This yields 
\[
\frac{1}{\rho} \, 
q_t 
=
-2 \, v \, v_x  
-
\frac{\rho_x }{\rho} \, v^2 +
\frac{1}{\rho} \, 
\left(
\sum_{i=1}^n  
- \rho^i \,  P_i''(\rho^i)\, \rho^i_x 
-
\sum_{i=1}^n  
\left(
\rho^i\, (u^i)^2
\right)_x 
-
\sum_{i=1}^n  
 \frac{1}{2}\theta \, \rho^i \,   v^i\, |v^i|
 \right)
.
\]
Due to  (\ref{vtgleichung})  this yields
\[
v_t
=
\frac{1}{\rho} \, q_t  + v \, v_x + v^2 \, \frac{\rho_x}{\rho}
\]
\[
=
-  v \, v_x +
\frac{1}{\rho} \, 
\left(
\sum_{i=1}^n  
- \rho^i \,  P_i''(\rho^i)\, \rho^i_x 
-
\sum_{i=1}^n  
\left(
\rho^i\, (u^i)^2
\right)_x 
-
\sum_{i=1}^n  
 \frac{1}{2}\theta \, \rho^i \,   v^i\, |v^i|
 \right)
 .
\]
Thus we have
\[
v_t + v \, v_x +
\left(
\sum_{i=1}^n  
 \frac{ \rho^i }{\rho}  \,  P_i''(\rho^i)\, \rho^i_x 
\right)
=
\frac{1}{\rho} \, 
\left( 
-
\sum_{i=1}^n  
\left(
\rho^i\, (u^i)^2
\right)_x 
-
\sum_{i=1}^n  
 \frac{1}{2}\theta \, \rho^i \,   v^i\, |v^i|
 \right).
\]
Thus we have also shown the second equation   in  (\ref{hamiltonianmixturn}).
\end{proof}
%WEG !!!!!!!!!!!!!
%In the case of ideal gases, where $\rho^i  \,  P_i''(\rho^i) = a_i $, this yields
%\[
%v_t +  \left( \frac{1}{2} v^2 \right)_x+  
%\frac{1}{\rho} \sum_{i=1}^n a^i  \, \rho^i_x =  ...
%\]

%Using (\ref{mischfluss}) this yields the second equation in  (\ref{hamiltonianmixture}).

%%%%%%%%%%%%%%%%%%%%%%%%%%%%%%%%%%%%%%%%%%%%%%%%%%%%%%%%%%%%%%%%%%%%%%%%%%%%%%%%%%%%%%%%%%%%%%%%%%%%%%%%%%%%%%%%%%%%%%%%%%%%%%

%%%%%%%%%%%%%%%%%%%%%%%%%%%%%%%%%%%%%%%%%%%%%%%%%%%%%%%%%%%%%%%%%%%%%%%%%%%%%%%%%%%%%%%%%%%%%%%%%%%%%%%%%%%%%%%%%%%%%%%%%%%%%%

\section{The Lyapunov function}
\label{Lyapunovmix}

The proof of the exponential synchronization 
of the velocities for appropriate boundary data  is based on  a suitably chosen Lyapunov function.
Our choice of the Lyapunov function is
motivated by relative energy techniques that have
been used for example in \cite{Dafermos_book},  \cite{diperna}, \cite{herbertegger},  \cite{GiesselmannGugatKunkel}.
However, in contrast to \cite{Dafermos_book},  \cite{diperna}, \cite{herbertegger},  \cite{GiesselmannGugatKunkel}  
in our Lyapunov function the derivative of the energy functional does not appear.
For  $t\geq 0$,   consider the Lyapunov function
\[
{\cal L}^e(t) =  \int_0^{L}   
\left( \sum_{i=1}^n H^i(\rho^i , \, v^i) \right)  - H^{mix}(\rho^1, \ldots, \rho^n, \,  v)
\, dx.
\]

Then  we have the equation 
\begin{equation}
    \label{lerepresentationmix}
{\cal L}^e(t) = %\frac{1}{2}  \, 
\int_0^{L} 
 \sum_{i=1}^n 
\frac{1}{2} \,  \rho^i \, |  v^i -  v|^2 \, dx.
\end{equation}
This can be seen as follows.
For $i \in \{1,\ldots, n\}$,  let 
\begin{equation}
    \label{lambdaidefi}
\lambda_i = \frac{  \rho^i }{  \rho  } .
\end{equation}

For all $i \in \{1,\ldots, n\}$ we have
\[
H^i(\rho^i , \, v^i) = \frac{1}{2} \, \rho^i \, (v^i)^2  + P_i
=
\frac{1}{2} \rho^i \, \left( v + (v^i-v)  \right)^2 + P_i
\]
\[
=
\frac{1}{2} \rho^i \,  v^2 + \rho^i \, v \, (v^i - v ) +  \frac{1}{2} \rho^i  \, (v^i-v)^2 + P_i
\]
\[
=
\lambda_i 
 ( H^{mix}(\rho^1, \ldots, \rho^n  , \, v)  - \sum_{j=1}^n P_j )
+
\rho^i \, v \,
u^i
+
\frac{1}{2} \, 
\rho^i \, 
(u^i)^2
+
P_i
.
\]
Equation (\ref{rhoiuigleichung}) states that 
\[
\sum_{i=1}^n  \rho^i \, u^i =0.
\]
%Note that  we have
%\begin{equation}
%    \label{nullsumme}
%\sum_{i=1}^n  \rho^i \,  (v^i - v )  = \left(\sum_{i=1}^n  \rho^i \,  v^i \right)  -v\, \sum_{i=1}^n  \rho^i 
%= \rho \, v - v  \, \rho = 0.
%\end{equation}
Hence summation for $i \in \{1,\ldots,n\}$ yields (\ref{lerepresentationmix}). To be precise, 
%Due to convexity 
 we have
\[
\sum_{i=1}^n H^i(\rho^i , \, v^i) 
= 
H^{mix}(\rho^1, \ldots, \rho^n, \,  v)
+ \sum_{i=1}^n \frac{1}{2} \, \rho_i \,  (v_i - v)^2.
\]
Note that these arguments are similarly as in the derivation of  Jensen's inequality for strongly convex functions,
see  e.g. \cite{bakula2016converse}. Strongly convex functions are also discussed in
\cite{stroco} and \cite{vial}.
%On the converse Jensen inequality for strongly convex functions
%Milica Klaričić Bakula a
%, Kazimierz Nikodem 
%2026

%%%%%%%%%%%%%%%%%%%%%%%%%%%%%%%%%%%%%%%%%%%%%%%%%%%%%%%%%%%%%%%%%%%%%%%%%%%%%%%%%%%%%%%%%%%%%%%%%%%%%%%%%%%%%%%%%%%%%%%%%%%%%%

\section{Evolution of the Lyapunov function}
\label{evolutionmix}
In the next step, we study how  the Lypunov function evolves with time. 
%This is done by looking at the derivative.
For the time-derivative of   $ {\cal L}^e  $  we have 
\[
\frac{d}{dt}
{\cal L}^e(t)
= \int_0^{L}   
 \sum_{i=1}^n 
 \left[
 H^i_{\rho^i}(\rho^i , \, v^i)\, \rho^i_t
 +
  H^i_{v^i}(\rho^i , \, v^i)\, v^i_t 
  \right]
  \]
  \[
-
 \sum_{i=1}^n  
 \left[
 H^{i}_{\rho^i}(\rho^i , \, v)\, \rho^i_t
 \right] 
  - 
 H^{mix}_v(\rho^1, \ldots, \rho^n, \,  v) \, v_t
\, dx.
\]
Due to (\ref{hamiltonianmix}),  (\ref{hamiltonianmixturn})
and (\ref{varderi}), (\ref{vableitung})  this yields
\[
\frac{d}{dt}
{\cal L}^e(t)
=
 \int_0^{L}   
  \sum_{i=1}^n 
   \left[
  - H^i_{\rho^i}(\rho_i , \, v_i)\, \partial_x  H^i_{v^i} (\rho^i, \, v^i)
 -
  H^i_{v^i}(\rho^i , \, v^i)\, \partial_x  H^i_{\rho^i} (\rho^i, \, v^i)
   \right]
\]
\[
- 
\sum_{i=1}^n  
 \left[
\rho^i \, \frac{1}{2}\theta \,   (v^i)^2\, |v^i|
-
%(\beta + M   + \frac{1}{2} \theta \, \,{|{v}|} \  )
\bar\Omega \,  \rho^i   ( v^i - v)^2
 \right]
\]
%%\[
%%????????
%%+ \rho^i \, (v^i -v)   \,  P_i''(\rho^i)\, \rho^i_x 
%%\]
%%Definition  von $\bar G^i$ anpassen.
 \[
+
 \sum_{i=1}^n  
  \left[
 \left( 
 \frac{1}{2} v^2 + P_i'(\rho^i)
 \right)
 \,\partial_x  (\rho^i v^i)
  \right]
\]
  \[
 +
\rho \, v  \, 
\left[
v \, v_x 
+
\left(
\sum_{i=1}^n  
 \frac{ \rho^i }{\rho}  \,  P_i''(\rho^i)\, \rho^i_x 
\right)
+
\frac{1}{\rho} \, 
\left( 
\sum_{i=1}^n  
\left(
\rho^i\, (u^i)^2
\right)_x 
+
\sum_{i=1}^n  
 \frac{1}{2}\theta \, \rho^i \,   v^i\, |v^i|
 \right)
 \right]
 \, dx.
 \]
 %%%%%%%%%%%%%%%%%%%%%%%%%%%%%%%%%%%%%%%%%%%%%%%%%%%%%%%%%%%%%%%%%%%%%%%%%%%%%%%%%%%%%%%%%%%%%%%%%%
Thus  using the product rule and  (\ref{rhovgleichung})   we obtain 
\[
\frac{d}{dt}
{\cal L}^e(t)
=
 \int_0^{L}   
  \sum_{i=1}^n 
   - \partial_x  \left( H^i_{\rho^i}(\rho_i , \, v_i)\,  H^i_{v^i} (\rho^i, \, v^i) \right)
\]
\[
- 
\sum_{i=1}^n  
 \left[
\rho^i \, \frac{1}{2}\theta \,   (v^i)^2\, |v^i|
-
%(\beta + M   + \frac{1}{2} \theta \, \,{|{v}|} \  ) 
\bar\Omega 
\,  \rho^i   ( v^i - v)^2
 \right]
\]

 \[
+
 \sum_{i=1}^n  
\partial_x
\left(
 \left( 
 \frac{1}{2} v^2 + P_i'(\rho^i)
 \right)
  (\rho^i v^i)
  \right) 
   -  \sum_{i=1}^n  
\left(  
u^i \,   \rho^i  \,  P_i''(\rho^i)\, \rho^i_x 
\right)
\]
  \[
 +
 v  \, 
\left[
\sum_{i=1}^n  
\left(
\rho^i\, (u^i)^2
\right)_x 
+
\sum_{i=1}^n  
 \frac{1}{2} \, \theta \, \rho^i \,   v^i\, |v^i|
 \right]
 \, dx.
 \]
 %%%%%%%%%%%%%%%%%%%%%%%%%%%%%%%%%%%%%%%%%%%%%%%%%%%%%%%%%%%%%%%%%%%%%%%%%%%%%%%%%%%%%%%%%%%%%%%%%%%%%%%%
Define the function 
\[
B(t)=
\int_0^{L}   
  \sum_{i=1}^n 
   - \partial_x  \left( H^i_{\rho^i}(\rho_i , \, v_i)\,  H^i_{v^i} (\rho^i, \, v^i) \right)
+
 \sum_{i=1}^n  
\partial_x
\left(
 \left( 
 \frac{1}{2} v^2 + P_i'(\rho^i)
 \right)
  (\rho^i v^i)
  \right) \,dx
\]
\begin{equation*}
%    \label{Bdefinition}
=\int_0^{L}   
 \sum_{i=1}^n  
\partial_x
\left(
\left( 
 \frac{1}{2} v^2  -    \frac{1}{2} (v^i)^2 
 \right)(\rho^i \, v^i) 
 \right)
 \,dx.
\end{equation*}
\begin{equation}
    \label{Bdefinition}
=
 \sum_{i=1}^n  
\left(
\left( 
 \frac{1}{2} v^2  -    \frac{1}{2} (v^i)^2 
 \right)(\rho^i \, v^i) 
 \right)|_{x=0}^L
.
\end{equation}
%%%%%%%%%%%%%%%%%%%%%%%%%%%%%%%%%%%%%%%%%%%%%%%%%%%%%%%%%%%%%%%%%%%%%%%%%%%%%%%%%%%%%%%%%%%%%%%%%%%%%%%%%%%%%%%%%%%%
 Then using  integration by parts we obtain the equation 
 \[
\frac{d}{dt}
{\cal L}^e(t)
=
B(t) +  \int_0^{L}   
- 
\sum_{i=1}^n  
 \left[
\rho^i \, \frac{1}{2} \, \theta \,   (v^i)^2\, |v^i|
-
\bar\Omega 
%(\beta + M   + \frac{1}{2} \, \theta \, \,{|{v}|} \  ) 
\,  \rho^i   ( v^i - v)^2
 \right]
\]

 \[
   -  \sum_{i=1}^n   
u^i \,   \rho^i  \,  P_i''(\rho^i)\, \rho^i_x 
 +
v \sum_{i=1}^n  
 \frac{1}{2}\theta \, \rho^i \,   v^i\, |v^i|
-  v_x \, \sum_{i=1}^n  
 \left(
\rho^i\, (u^i)^2
\right)
 \, dx
 +
\left[v \, \sum_{i=1}^n  
 \left(
\rho^i\, (u^i)^2
\right)
\right]_{x=0}^L.
 \]
Define 
\begin{equation}
\label{Bhutdefinition}
\hat B(t):= 
B(t) 
+
\left[v \, \sum\limits_{i=1}^n  
 \left(
\rho^i\, (u^i)^2
\right)
\right]_{x=0}^L.
\end{equation}
Then we have the representation 
\[
\frac{d}{dt}
{\cal L}^e(t)
=
\hat B(t) +  \int_0^{L}   
 \frac{1}{2}\theta \,
 \sum_{i=1}^n  
 \left[
 (v - v^i)
 \rho^i \,   v^i\, |v^i|
 \right]
-
\sum_{i=1}^n  
 \left[
\bar\Omega
%(\beta + M   + \frac{1}{2} \, \theta \, \,{|{v}|} \  )
\,  \rho^i   ( v^i - v)^2
 \right]
\]
 \[
   -  \sum_{i=1}^n   
   \left[
u^i \,   \rho^i  \,  P_i''(\rho^i)\, \rho^i_x 
\right]
-  v_x \, \sum_{i=1}^n  
 \left[
\rho^i\, (u^i)^2
\right]
 \, dx
 .
 \]

%%%%%%%%%%%%%%%%%%%%%%%%%%%%%%%%%%%%%%%%%%%%%%%%%%%%%%%%%%%%%%%%%%%%%%%%%%%%%%%%%%%%%%%%%%%%%%%%%%%%%%%%%%%%%%%%%%

\subsection{An estimate for the friction terms}
Now we consider  the terms  in our  representation of 
$
\frac{d}{dt}
{\cal L}^e(t)
$
that are generated by the terms  that
model the  friction at the   
boundary of the  pipeline.
As long as the sign of the velocity does not change,
these terms can be considered as quadratic. 
%where the friction parameter $\theta$
We define
\[
U(t) = 
%\sum_{i=1}^n  
% \left[
%\rho^i \,   (v^i)^2\, |v^i|
%-
% \rho^i \,
%v \,  v^i\, |v^i|
% \right]
% =
 \sum_{i=1}^n  
\rho^i \,(v^i - v)    \, v^i\, |v^i|
.
\]
Due to (\ref{rhovgleichung})   we have
\[
\sum_{i=1}^n   \rho^i \, v \,(v^i - v) =0,
\]
and thus
\begin{equation}
    \label{hilfsgleichung4}
\sum_{i=1}^n   \rho^i \, v^i \,(v^i - v) 
=
\sum_{i=1}^n   \rho^i \,\,(v^i - v)^2 
.
\end{equation}
%Moreover, we obtain
%\[
%\sum_{i=1}^n   \rho^i \, (v^i)^2 \,(v^i - v) 
%=
%\sum_{i=1}^n   \rho^i \,\,(v^i - v)^3    + 2 \sum_{i=1}^n   \rho^i \,v \,  (v^i - v)^2
%=
%\sum_{i=1}^n   \rho^i \,( v + v^i) \,(v^i - v)^2  
%.
%\]

%%%%%%%%%%%%%%%%%%%%%%%%%%%%%%%%%%%%%%%%%%%%%%%%%%%%%%%%%%%%%%%%%%%%%%%%%%%%%%%%
Assume that 
$ \sigma = sign( v^i )$ is independent of $i$ 
(this means that the flow of all components has the same orientation,
for example $v^i>0$ for all $i \in \{1,\ldots,n\}$). 
Then  with (\ref{hilfsgleichung4}) we obtain 
\[
U(t)
=\sigma 
\sum_{i=1}^n  
\rho^i \,(v^i - v)   \, v^i\, ( v + u^i)
=
\sigma 
\sum_{i=1}^n  
\rho^i \,(v^i - v)   \, v^i\, v
+
\sigma 
\sum_{i=1}^n  
\rho^i \,(v^i - v)  \, v^i \, u^i
\]
\[
=
\sigma 
\sum_{i=1}^n  
\rho^i 
\,
v\,(v^i - v)^2
+
\sigma 
\sum_{i=1}^n  
\rho^i \, v^i\,(v^i - v)^2
=
\sigma 
\sum_{i=1}^n  
\rho^i 
\,
(v + v^i)
\,(v^i - v)^2.
\]

This yields the inequality
\[
|U(t)|\leq 
\sum_{i=1}^n  
\rho^i 
\,
|v + v^i|
\,(v^i - v)^2.
\]

%%
%WEG:
%(ERSETZE DAHER in $\bar G^i$ die Zahl  $\frac{1}{2} \theta $ durch $\theta$!!!!!!!!!!!
%Oder $|v|$ durch $|v + v^i|$!!!?? Geht nicht, da dann dass Aufsummieren zu Null kaputt geht!!)

%Define the auxiliary function $h(x) = x \, |x|. $
%Then $f$ is continuously differentiable and $h'(x) = 2\,|x|$.
%Moreover, for $x\not=0$, we have $h''(x) = 2\, sign(x)$.
%Thus for $v\not=0$, there exist a number $\eta^i \in [-1, \, 1]$ such that we have
%\[
%v^i |v^i| = h(v^i)
%= h(v + u^i)
%=  h(v) + h'(v) \,u^i  +  \eta^i \, (u^i)^2.
%\]
 
%%%%%%%%%%%%%%%%%%%%%%%%%%%%%%%%%%%%%%%%%%%%%%%%%%%%%%%%%%%%%%%%%%%%%%%%%%%%%%%%%%%%%%%%%%%%%%%%%%%%%%%%%%%%%%%%%%

\subsection{A differential  inequality}
Now we derive a differential inequality for our Lyapunov function.

{%\color{blue}
Let  constants $M>0$, $N>0$
and $ \hat\varepsilon >0$ 
be given 
%that are sufficiently large 
such that 
%Assume  that the constants $M$ and $N$ are sufficiently large such that 
for all $i \in \{1,\ldots,n\}$
we have the inequalities }
%almost everywhere 
%the velocity bound 
 \begin{equation}
     \label{Massv1}
\frac{\theta }{2} \,  |v^i| +   | \partial_x v| \leq M,
\end{equation}
\begin{equation}
\label{Massv1n}
 \frac{1}{2} \theta \,{|{v}|} \leq  N 
\end{equation}
and
 \begin{equation}
     \label{Massv2}
%\sum\limits_{i=1}^n 
\max_{i \in \{1, \ldots, n\}}
\rho^i \,   | P_i''(\rho^i)\, \rho^i_x |^2    
\leq  \frac{\hat\varepsilon^2}{2\, L}.
\end{equation}

Assume that  
\[  \beta :=  \bar \Omega - M - N >0  .\]
We have 
\[\bar \Omega = M + N + \beta.\]
Define
\begin{equation}
    \label{Ilerepresentationmix}
 I(t) = 
 \sum\limits_{i=1}^n 
\int\limits_0^{L}     
 %  \left[
(v^i-v) \,   \rho^i  \,  P_i''(\rho^i)\, \rho^i_x 
\,dx
 %\frac{1}{2}  \, 
%\int_0^{L} 
% \sum_{i=1}^n 
%\frac{1}{2} \,  \rho^i \, |  v^i -  v|^2 \, dx
.
\end{equation}

Then with $\hat B(t)$ as defined in (\ref{Bhutdefinition}) we  obtain the inequality 

$
\frac{d}{dt}
{\cal L}^e(t)
\leq 
\hat B(t) 
- 
\sum\limits_{i=1}^n 
\int\limits_0^{L}     
 %  \left[
(v^i-v) \,   \rho^i  \,  P_i''(\rho^i)\, \rho^i_x 
+
\beta  \,   \rho^i \, |v^i - v|^2 
%\right]
\, dx
$
\[=
\hat B(t) 
- 
I(t) 
-
2\, 
\beta
\,
{\cal L}^e(t)
\]
with $I(t)$ 
as defined in
(\ref{Ilerepresentationmix}).
%Thus we have
%\[
%\frac{d}{dt}
%{\cal L}^e(t)
%\leq 
%\hat B(t) 
%+
% \hat \varepsilon
% \,
% \left(
%\int\limits_0^{L}     
%\sum\limits_{i=1}^n 
%(v^i-v)^2
%\, dx 
%\right)^{1/2}
%-
%\beta 
%\int\limits_0^{L}     
%\sum\limits_{i=1}^n  
%\rho^i \,
%|v^i - v|^2 
%\, dx
%.
%\]
Hence  we have the differential inequalities
\begin{equation}
\label{iiungleichungI}
\frac{d}{dt}
{\cal L}^e(t)
\leq 
\hat B(t) - I(t)  -
2\, 
\beta
\,
{\cal L}^e(t).
\end{equation}
Since
\begin{equation}
    \label{Ileestimate}
 I(t) = 
 \sum\limits_{i=1}^n 
\int\limits_0^{L}     
 %  \left[
 \sqrt{\rho^i} (v^i-v) \,   \sqrt{\rho^i } \,  P_i''(\rho^i)\, \rho^i_x 
\,dx
 %\frac{1}{2}  \, 
%\int_0^{L} 
% \sum_{i=1}^n 
%\frac{1}{2} \,  \rho^i \, |  v^i -  v|^2 \, dx
\leq
\frac{\hat \varepsilon }{ \sqrt{2 }  } \, 
 \left(
\int\limits_0^{L}     
\sum\limits_{i=1}^n 
\rho^i \,
(v^i-v)^2
\, dx 
\right)^{1/2}
\end{equation}
inequality  (\ref{iiungleichungI}) implies  
\begin{equation}
\label{iiungleichung}
\frac{d}{dt}
{\cal L}^e(t)
\leq 
\hat B(t) +
 \hat \varepsilon
 \,
 \left(
{\cal L}^e(t) 
%\int\limits_0^{L}     
%\sum\limits_{i=1}^n 
%(v^i-v)^2
%\, dx 
\right)^{1/2}
-
2\, 
\beta
\,
{\cal L}^e(t)
.
\end{equation}

If $\hat B(t) \leq 0$, this implies
\begin{equation}
\label{zentralungleichung}
\frac{d}{dt}
{\cal L}^e(t)
\leq 
 \hat \varepsilon
 \,
 \left(
{\cal L}^e(t)
%\int\limits_0^{L}     
%\sum\limits_{i=1}^n 
%(v^i-v)^2
%\, dx 
\right)^{1/2}
-
2\, 
\beta
\,
{\cal L}^e(t)
.
\end{equation}
%If we assume that (\ref{rhoischranke}) holds,  we obtain
%\begin{equation}
%\label{zentralungleichung}
%\frac{d}{dt}
%{\cal L}^e(t)
%\leq 
% \frac{\hat \varepsilon}
% {M_0} \,
%\sqrt{
%{\cal L}^e(t)
%}
%-
%2\, 
%\beta
%\,
%{\cal L}^e(t)
%.
%\end{equation}
\begin{remark}
It is important to note that with (\ref{zentralungleichung}),
the initial condition ${\cal L}^e(0)=0$ does not imply
${\cal L}^e(t)=0$ for $t\geq 0$.
This can be seen as follows:
Consider the differential equation
\[y' =  
\alpha \, \sqrt{y} - 2 \beta \, y.
 \]
 with $\alpha  =   \frac{\hat \varepsilon}
 {M_0}   $.
 The right-hand side is not Lipschitz continuous as
 a function of $y$ at zero.
 With the initial condition $y(0)=0$,
 the differential equation has several solutions,
for example  $y(t)=0$ and
another solution is
\[y(t) =  \frac{\alpha^2 }{4\, \beta^2 } \, \left( 1 -  \exp(  - \beta \, t)    \right)^2.\]
For the latter we have $y(0)=0$, $y'(0)=0$ and $y(t)>0$ for $t>0$.
More generally, for initial states $y(0)>  \frac{\alpha^2 }{4\, \beta^2 }  $
the differential equation has the solution
\[y(t) = 
\frac{\alpha^2 }{4\, \beta^2 } \, \left( 1  +  {C_0} \exp(  - \beta \, t)    \right)^2
\]
with $C_0 > -1$.
Thus we see that for all initial states $y(0)>0$,
for the state we have
\[\lim_{t \rightarrow \infty} y(t) =  \frac{\alpha^2 }{4\, \beta^2 }  \]
and the convergence is exponentially fast with the rate $\beta$.

\end{remark}

\subsection{Negativity of the terms that are  evaluated at the boundary points}
Due to 
(\ref{Bhutdefinition}) and  (\ref{Bdefinition}) 
we have
\begin{equation}
\label{Bhutdefinitionfolge}
\hat B(t)= 
\sum_{i=1}^n  
\left[
\left( 
 \frac{1}{2} v^2  -    \frac{1}{2} (v^i)^2 
 \right)(\rho^i \, v^i) 
+
v \, 
%\sum\limits_{i=1}^n  
 \left(
\rho^i\, (u^i)^2
\right)
\right]_{x=0}^L
\end{equation}
\[
=
\sum_{i=1}^n  
\left[
\left( 
 \frac{1}{2} (v + v^i)( v- v^i)  
 \right)(\rho^i \, v^i) 
+
v \, 
%\sum\limits_{i=1}^n  
 \left(
\rho^i\, (v - v^i )^2
\right)
\right]_{x=0}^L
\]
\[
=
\sum_{i=1}^n  
\rho^i
\left[
\left( 
 \frac{1}{2} (v + v^i)
\, v^i 
+
v \,  (v - v^i )
\right)
 (v - v^i )
\right]_{x=0}^L
\]
\[
=
\sum_{i=1}^n  
\rho^i
\left[
\left( 
 v^2 -   \frac{1}{2} v\,  v^i
+ \frac{1}{2}    (v^i )^2
\right)
 (v - v^i )
\right]_{x=0}^L.
\]
{%\color{blue}
The  definition of the barycentric velocity in (\ref{vdefinition}) 
implies that $\sum_i \rho^i \, (v^i - v) = 0$.}
Hence we have
%this yields
\[
\hat B(t) =
\sum_{i=1}^n  
\rho^i
\left[
\left( 
-   \frac{1}{2} v\,  v^i
+ \frac{1}{2}    (v^i )^2
\right)
 (v - v^i )
\right]_{x=0}^L
\]
\[
=
\sum_{i=1}^n  
\rho^i
 \frac{1}{2}   v^i
\left[
\left( 
-   v
+   v^i )
\right)
 (v - v^i )
\right]_{x=0}^L
\]
\[
=
-
 \frac{1}{2} 
\sum_{i=1}^n  
\rho^i
\left[
  v^i
\,
 (v - v^i )^2
\right]_{x=0}^L .
%<0 ???????????????????????
\]
%Since $\sum_i \rho^i v^i = \rho \, v$, this yields
%\[
%\hat B(t) 
%=
%\left[
%\rho\, v^3 -   \rho \, \frac{3}{2} v^3
%+
%\frac{1}{2}
%\sum_{i=1}^n  
%\rho^i
%\, (v^i )^3
%\right]_{x=0}^L
%\]
%\[
%=
%\left[
%\frac{1}{2}
%\sum_{i=1}^n  
%\rho^i
%\,
%\left( (v^i )^3 - v^3 \right)
%\right]_{x=0}^L.
%\]

Define
\begin{eqnarray*}
B_L(t) & = &
%\sum_{i=1}^n  
%\left[
%\rho^i
%\left( 
% v^3 -   \frac{3}{2} v^2\,  v^i
%+ \frac{1}{2}    (v^i )^3
%\right)
%\right](t, \, L)
%=
-
 \frac{1}{2} 
\sum_{i=1}^n  
\rho^i
  v^i
\,
 (v - v^i )^2
(t,\, L)
%
%\frac{1}{2}
%\sum_{i=1}^n  
%\rho^i
%\,
%\left( (v^i )^3 - v^3 \right)(t,\, L)
,
\\
B_0(t) & = &
%\sum_{i=1}^n  
%\left[
%\rho^i
%\left( 
% v^3 -   \frac{3}{2} v^2\,  v^i
%+ \frac{1}{2}    (v^i )^3
%\right)
%\right](t, \, 0).
%=
-
 \frac{1}{2} 
\sum_{i=1}^n  
\rho^i
  v^i
\,
 (v - v^i )^2
(t,\, 0).
\end{eqnarray*}
Then we have
\[\hat B(t) = B_L(t) - B_0(t).\]
%We can represent $B_0(t)$ in the form 
%\[
%B_0(t) 
%=
%\frac{\rho}{2} \,
%\left(
%\sum_{i=1}^n  
%\lambda^i
%\,
% (v^i )^3 - 
% \left(
%\sum_{j=1}^n  
%\lambda_j v^j
%\right)^3(t, \, 0)
%\right) 
%\]
%where the $\lambda_i$ are as in (\ref{lambdaidefi}).
%Since the  function  $w\mapsto w^3$ is convex on $[0, \infty)$,
%if $v^i \geq 0$ for all $i \in \{1,\ldots, n\}$ 
If $v^i \geq 0$ for all $i \in \{1, \ldots,n\}$, 
this yields the inequalities
\[
B_L(t) \leq 0, \; 
 B_0(t)  \leq 0.
\]
%Since the maps $w\mapsto w^3$ is strictly convex on $[0, \infty)$,
If $\sum_{i=1}^n (u^i)^2(t, \, L)  >0$ 
%and $\lambda_i \in (0, 1)$ 
%for all $i \in \{1,\ldots,n\}$
the strict inequality 
\[
 B_L(t)  < 0
\]
holds.
If $v^i(t,\, 0) = v(t, 0) $ for all $i \in \{1,\ldots,n\}$ we have $B_0(t)=0$.
This is the case with the  boundary condition  (\ref{5.10}).

%Similarly, we have $B_L(t) \geq 0$ and with the control law
%(\ref{5.10})  
%%real number
%we have
%$B_L(t)=0$.

%%%%%%%%%%%%%%%%%%%%%%%%%%%%%%%%%%%%%%%%%%%%%%%%%%%%%%%%%%%%%%%%%%%%%%%%%%%%%%%%%%%%%%%%%%%%%%%%%%%%%%%%%%

\section{The synchronization results}
\label{sec:sync}

 Theorem \ref{thm7.1} 
 contains a sufficient  condition for exponential synchronization
of the velocities with suitable boundary data. 
In the following remark we discuss the assumptions:
\begin{remark}
{%\color{blue}
Theorem \ref{semiglobalphysical} implies
that in a $C^1$-neighbourhood with radius $ \varepsilon(T) $ of
a stationary solution of 
(\ref{massi}), 
(\ref{momenti}), 
(\ref{5.10}), 
(\ref{feedbackrho}) 
we can find classical solutions
of the corresponding initial boundary value problem
that 
satisfy the a-priori inequality (\ref{aprioribound2025}).

To be precise,  let us have a closer look at the stationary solutions.
For the stationary solutions, we use the notation
$\tilde \rho^i$, $\tilde v^i$.
Then for all $i \in \{1, \ldots, n\}$, the flow rate 
\[\tilde q^i = \tilde \rho^i(x) \, \tilde v^i(x)\]
is constant
for $x \in [0, \, L]$,
for $x_B=0$ 
we have $\tilde v^i(x_B) = \bar v$
and 
%at $L$ 
we have $\tilde \rho^i(L) = \bar \rho^i$.
We have 
$\tilde v^i = \frac{\tilde q^i}{\tilde \rho^i}$
and
$ v= \frac{\sum_i \tilde q^i}{\sum_i \tilde \rho^i}$.
The values of the functions $\tilde \rho^i(x)$ for
$x\in  [0, \, L)$ 
and $\tilde v^i(x)$ for $x \in (0, \, L]$
are
determined by the ordinary differential equations
\[
\left(  p^i + \frac{(\tilde q^i)^2}{{ \tilde  \rho}^i} \right)_x = 
- \frac{1}{2} \theta \, \frac{\tilde{q}^i \; |\tilde{q}^i|}{\tilde{\rho}^i}  -
\bar \Omega
%\left( \beta  + M
 %\frac{1}{2} \theta%{\rho} \,{|{v}|} 
%+ N   \right)
%\,  \rho^i \,  ( v^i - v) 
\left(\tilde q^i - 
  \frac{ \sum_j {\tilde q^j} }{ \sum_j {\tilde \rho^j}} \, 
\tilde \rho^i
\right)
\]
where the numbers $\tilde q^i$ are chosen such that
for all $i \in \{1,\ldots,n\}$ we have 
$
\tilde \rho^i(0) =\frac{ \bar v}{ \tilde q^i   }.
$

In particular, there are  the  
constant stationary states with zero velocity, that is 
 for all $i \in \{1, \ldots, n\}$
 we have $\tilde v^i(x)=0$
and the  constant densities $\tilde \rho^i(x)  = \bar \rho^i$
($x\in [0, L]$).

Let  constants $M>0$, $N>0$
and $ \hat\varepsilon >0$ be given 
such that  $\beta = \bar\Omega - N - M >0$.
The a-priori inequality (\ref{aprioribound2025})
implies that if  the number $\varepsilon(T)>0$
is chosen sufficiently small,
the classical solutions
in the  $C^1$-neighbourhood with radius $ \varepsilon(T) $
of
a stationary state 
$\tilde \rho^i$, $\tilde v^i >0$  
satisfy 
(\ref{Massv1}),  (\ref{Massv1n})  and (\ref{Massv2}).
%and  $\beta = \bar\Omega - N - M >0$.

In a nutshell, this means that
the assumptions of Theorem \ref{thm7.1}
that is stated below 
hold 
for 
$C^1$-compatible  initial data and   boundary data
in a sufficiently small  $C^1$-neighbourhood 
of  stationary states for the system
with positive velocities. 
}
\end{remark}

\begin{theorem}
\label{thm7.1}
%\begin{itemize}
Let $T>0$ be given. 
Assume that the 
%coupled 
system  that consists of 
(\ref{massi}), 
(\ref{momenti}), 
(\ref{initial}), 
(\ref{5.10}), 
(\ref{feedbackrho}) 
for all $i\in \{1,\ldots,n\}$
has a classical solution on $[0, \, T]$ such that 
(\ref{Massv1}),  (\ref{Massv1n})  and (\ref{Massv2}) hold
and  $\beta = \bar\Omega - N - M >0$.
Assume that 
$  v^i \geq 0 $
%and $\rho^i >0$ 
for all  $i \in \{1,\ldots,n\}$.

           (i)
%Assume that in addition 
% (\ref{rhoischranke}) holds.
% and  for all $t\geq 0$ we have  $\hat B(t) \leq 0$
% that is $B_L(t) \leq B_0(t)$. 
% This is the case for example if $B_0(t)=0$,  that is  with the boundary conditions  (\ref{5.10}). 
The values of  
    $
    {\cal L}^e(t)
    $
    (as in (\ref{lerepresentationmix}) ) 
    decay exponentially with time 
    %fast with the rate $2\beta$, 
    as long as   ${\cal L}^e(t)$ 
    has values that are 
    outside of  the interval
    \[
    I_{syn} =
    \left[0, \, \frac{ \hat \epsilon^2 }{ 4 \,   \beta^2 }   
    %\frac{ \hat \epsilon^2}{M_0^2}  
    \right]
    \]
    %is of the order $ \frac{1}{\beta^2}$ 
    in the sense that for all $t \in (0,T) $  we have
    \begin{equation}
        \label{iii}
{{\cal L}^e(t) } 
\leq\left( \frac{ \hat \epsilon}{2\,  \beta}
+ 
\exp(- \beta\, t) \,
\left(  \sqrt{{\cal L}^e(0) } - \frac{ \hat \epsilon}{2\,  \beta} \right) \right)^2
 .
         \end{equation}

    (ii)  
    If there is $t_\ast \in [0, T]$ such that for all $t \in  [t_\ast,  T] $ we have  

   \begin{equation}
    \label{decayassumptionI}
    \hat B(t) 
    \leq  I(t),
    \; i.e. \;\;
  I(t) \geq B_L(t) 
\end{equation}
with $I(t)$ as defined in
(\ref{Ilerepresentationmix}), 
%or
%    \begin{equation}
%    \label{decayassumption}
%    \hat B(t) 
%    \leq -
% \hat \varepsilon
% \left(
%\int\limits_0^{L}     
%\sum\limits_{i=1}^n 
%(v^i-v)^2(t, \, x)
%\, dx 
%\right)^{1/2}
%=
%-
% \hat \varepsilon\, \sqrt{ 2 \,   {\cal L}^e(t) }
%\end{equation}
then 
$
    {\cal L}^e(t)
    $
    decays exponentially  fast  on $[t_\ast, \, T]$  with the rate $2\beta$,
    that is
    for all $ t\in [t_\ast,  T] $ we have
    \begin{equation}
\label{expdecay}
{\cal L}^e(t)
\leq 
\exp(- 2 \beta (t - t_\ast)  )\,{\cal L}^e(t_\ast). 
 \end{equation}
%(Here it is not assumed that  (\ref{rhoischranke}) holds).
 %\item

%\[
%\sqrt{{\cal L}^e(t) } - \frac{ \hat \epsilon}{2\, M_0\, \beta}
%\leq
%\exp(- \beta\, t) \,
%\left(  \sqrt{{\cal L}^e(0 } - \frac{ \hat \epsilon}{2\, M_0\, \beta} \right).
%  \]       
%\[
%    {\cal L}^e(t) 
%      \leq
%      \frac{ \hat \epsilon^2}{M_0^2} \,  \frac{1}{ ( 1 + \lambda)^2 \,   \beta^2 } + 
%    \max \left\{  {\cal L}^e(0) -   \frac{ {\hat \epsilon}^2}{M_0^2} \,  \frac{1}{ ( 1 +  \lambda)^2\, \beta^2 }, \, 0  \right\} \, \exp( - 
%    ( 1 - \lambda)  \beta \, t). 
%  \] 

%\linebreak 
    %\item

%\linebreak 
% (iii) 
% If $I(0)$ is sufficiently large (for example, if $I(0)\geq 0 $),  
% there exist boundary conditions such that 
% (\ref{decayassumptionI}) holds for all $t\in [0, \, T]$ 

% (for example if
%$  \underline C \, (\bar v(t) )^3\, u(t)^2 \geq - I (t)   $, 
%with the notation as in Lemma \ref{randbed}) 
 
% and
% hence the velocities of the different components 
% synchronize   exponentially  fast on $[0, \, T]$. 

%   (iv)
%  If there exists $t_0>0$ such that  $u^i(t) = 0$ for $i\in \{1,\ldots,n\}$,
%  then 
%  $\frac{d}{dt} {\cal L}^e(t_0) = 0$. 
  %$v^i(t) = v(t) $ for all $t\geq 0$.
    
%\end{itemize}
    
\end{theorem}

%\begin{remark}
%Condition (\ref{decayassumptionI}) holds if ???
%\end{remark}

\begin{remark}
 The upper limit of the synchronization interval  $I_{syn}$ is
\[   S_0  = \frac{ \hat \epsilon^2 }{ 4 \,   \beta^2 } 
%\,  \frac{ \hat \epsilon^2}{M_0^2}
.  \]
The denominator  is  of the order $\frac{1}{\beta^2}$.
So the size of the coupling constant 
$\bar \Omega$ dominates both the length of the interval  $I_{syn}$ 
and the decay rate. 
%The denominator also contains $M_0>0$, which is the lower bound 
%for the densities $\rho^i$.
%So the analysis in 
%\emph{(ii)} does not cover the case that certain components vanish.  
The numerator $\hat \varepsilon^2 $
can be chosen small if the derivatves $\rho^i_x$
(that is the pressure fluctuations) 
have sufficiently small absolute value, see the definition
of  $ \hat \varepsilon $  in  (\ref{Massv2}).

The result shows that the velocities will synchronize in
the $L^2$-sense until the
%error is 
synchronization error is less than or equal to $S_0$.
This means that 
the model error that is generated by replacing 
the different velocities $v^i$ by the average velocity $v$  
is bounded by $S_0$. This  implies 
a small model error if the parameter $ \bar\Omega $
of the velocity coupling is sufficiently large.

In fact, the simplified initial boundary value problem 
is
\begin{equation}
 \label{hamiltoniansimp}
\left\{
\begin{array}{l}
v(0, \, x) = v_0,
\;
\rho^i(0, \, x ) = \rho^i_0, \; x\in [0, \, L], \; i\in \{1,\ldots,n\}
\\
v(t, \, x_B) = \bar v(t),
\\
\rho^i(t, \, L)   = \bar \rho^i(t) ,\;   i \in \{1,\ldots,n\},
\\
 \rho^i_t +  (\rho^i \, v)_x = 0,\; i\in \{1,\ldots,n\},
 \\
 v_t + 
v \, v_x +
\left(
\sum_{i=1}^n  
 \frac{ \rho^i }{\rho}  \,  P_i''(\rho^i)\, \rho^i_x 
\right)
=
-
 \frac{1}{2}\theta \,  v\, |v|
.
 \end{array}
 \right.
\end{equation}
This is a drift flux model 
similar as 
for example in \cite{driftfluxbanda},  
\cite{driftflux}.

%consists of the initial conditions
%\begin{equation}
%\label{initialsimp}
%v(0, \, x) = v_0,
%\;
%\rho^i(0, \, x ) = \rho^i_0, \; x\in [0, \, L], \; i\in \{1,\ldots,n\}
%\end{equation}
%$\rho^i(t, \, L)   = \bar \rho^i(t) ,\;   i \in \{1,\ldots,n\}$
%(\ref{feedbackrho}),
%\begin{equation}
%\label{5.10velo}
%v(t, \, x_B) = \bar v(t),
%\end{equation}
%\begin{equation}
%    \label{massisimp}
% \rho^i_t +  (\rho^i \, v)_x = 0,
%\end{equation}
%\begin{equation}
%\label{vsimp}
%v_t + 
%v \, v_x +
%\left(
%\sum_{i=1}^n  
% \frac{ \rho^i }{\rho}  \,  P_i''(\rho^i)\, \rho^i_x 
%\right)
%=
%-
% \frac{1}{2}\theta \,  v\, |v|
%%\partial_x H^{mix}_\rho 
%\end{equation}

%for all $i\in \{1,\ldots,n\}$

%(\ref{massi}), 
%(\ref{momenti}), 
%(\ref{initial}), 

\end{remark}

\textbf{Proof of (i):}
 Consider the auxiliary function 
\begin{equation}
\label{Hdefinitionneu}
H(t) = \exp(  \beta \, t) \, \left(  \sqrt{{\cal L}^e(t)}  - \frac{ \hat \epsilon}{2\,  \beta}\right).
\end{equation}

Then  as long as  
$  {\cal L}^e(t)  >  0 $
%$ \frac{1}{(1 + \lambda)^2 }  \,  \frac{ {\hat \epsilon}^2}{M_0^2} \,  \frac{1}{  \beta^2 }  $
we have 
\[H'(t) =  \beta \, H(t) +  \exp(\beta \, t) \,  \frac{1}{ 2\,  \sqrt{{\cal L}^e(t)}   } \, \frac{d}{dt} \, {\cal L}^e(t). \]
Thus (\ref{zentralungleichung}) implies 
\[
H'(t)
\leq 
 \beta \, H(t) 
+  \exp(  \beta \, t) \,  \frac{1}{ 2\,  \sqrt{{\cal L}^e(t)}   }
\left(
- 2 \, \beta \, {\cal L}^e(t) +
{\hat \varepsilon}
 \,
\sqrt{
{\cal L}^e(t)
}
\right)
\]
\[
=
\beta \, H(t)  - \beta \, \exp( \beta \, t) \,  \sqrt{{\cal L}^e(t)}
+ \exp( \beta \, t) \, 
\frac{\hat \varepsilon}
 {2}
\]
\[
= 
 \beta \, H(t)
- \beta \,   \exp(\beta \, t) \,  
\left[
\sqrt{{\cal L}^e(t)} - \frac{\hat \varepsilon}
 {2 \, \beta }
\right]
= 0
\]
where the last equality fllows from (\ref{Hdefinitionneu}). 
Thus we have
$H(t) \leq H(0)$ and this implies
\[
\sqrt{{\cal L}^e(t) } - \frac{ \hat \epsilon}{2\, \beta}
\leq
\exp(- \beta\, t) \,
\left(  \sqrt{{\cal L}^e(0) } - \frac{ \hat \epsilon}{2\,  \beta} \right).
\]
This yields (\ref{iii}). 

%%%%%%%%%%%%%%%%%%%%%%%%%%%%%%%%%%%%%%%%%%%%%%%%%%%%%%%%%%%%%%%%%%%%%%%%%%%%%%%%%%%%%%%%%%%

\textbf{Proof of (ii)}: 
If 
(\ref{decayassumptionI})
holds,  inequality 
(\ref{iiungleichungI})
%(\ref{iiungleichung})
implies that for all $ t\in [t_\ast,  T] $ we have
\begin{equation}
\frac{d}{dt}
{\cal L}^e(t)
\leq 
-
2\, 
\beta
\,
{\cal L}^e(t).
\end{equation}
Then 
for all $ t\in [t_\ast,  T] $
 Gronwall's Lemma yields
 (\ref{expdecay})  
%\[
%{\cal L}^e(t)
%\leq 
%\exp(- 2 \beta (t - t_\ast)  )\,{\cal L}^e(t_\ast)  
%\]
and the assertion follows.
%with Gronwall's Lemma.

%\textbf{Proof of (iii)}:

%\textbf{Proof of (iv)}: If $u^i(t)=0$, we have
%$ {\cal L}^e(t) = 0 $,  $\hat B(t) = 0$ and $U(t)=0$.
%This yields $\frac{d}{dt} {\cal L}^e(t) = 0$.

%%%%%%%%%%%%%%%%%%%%%%%%%%%%%%%%%%%%%%%%%%%%%%%%%%%%%%%%%%%%%%%%%%%%%%%%%%%%%%%%%%%%%%%%%%%%%%%%%%%%%%%%%%%%%%%%%%%%

%%%%%%%%%%%%%%%%%%%%%%%%%%%%%%%%%%%%%%%%%%%%%%%%%%%%%%%%%%%%%%%%%%%%%%%%%%%%%%%%%%%%%%%%%%%%%%%%%%%%%%%%%%%%%%%%%%%%%%%%%%%%%%

%%%%%%%%%%%%%%%%%%%%%%%%%%%%%%%%%%%%%%%%%%%%%%%%%%%%%%%%%%%%%%%%%%%%%%%%%%%%%%%%%%%%%%%%

%%%%%%%%%%%%%%%%%%%%%%%%%%%%%%%%%%%%%%%%%%%%%%%%%%%%%%%%%%%%%%%%%%%%%%%%%%%%%%%%%%%%%%%%%%%%%%%%%%%%%%%%%%%%%%%%%%%%

%%%%%%%%%%%%%%%%%%%%%%%%%%%%%%%%%%%%%%%%%%%%%%%%%%%%%%%%%%%%%%%

%%%%%%%%%%%%%%%%%%%%%%%%%%%%%%%%%%%%%%%%%%%%%%%%%%%%%%%

%%%%%%%%%%%%%%%%%%%%%%%%%%%%%%%%%%%%%%%%%%%%%%%%%%%%%%%
%%%%%%%%%%%%%%%%%%%%%%%%%%%%%%%%%%%%%%%%%%%%

%%%%%%%%%%%%%%%%%%%%%%%%%%%%%%%%%%%%%%%%%%%%%%%%%%%%%%%%%%%%%%%%%%%%%%%%%%%%%%%%%%%%%%%%%%%%%%
%%%%%%%%%%%%%%%%%%%%%%%%%%%%%%%%%%%%%%%%%%%%%%%%%%%%%%%%%%%%%%%%%%%%%%%%%%%%%%%%%%%%%%%%%%%%%%

%%%%%%%%%%%%%%%%%%%%%%%%%%%%%%%%%%%%%%%%%%%%%%%%%%%%%%%%%%

%%%%%%%%%%%%%%%%%%%%%%%%%%%%%%%%%%%%%%%%

%%%%%%%%%%%%%%%%%%%%%%%%%%%%%%%%%%%%%%%%%%%%%%%%%%%%%%%%%%%%%%%%%%%%%%%%%%%%%%%%%%%%%%%%%%%%%%%%%%%%%%%%%%%%%%%%%%%%%%%%%

%%%%%%%%%%%%%%%%%%%%%%%%%%%%%%%%%%%%%%%%%%%%%%%%%%%%%%%%%%%%%%%%%%%%%%%%%%%%%%%%%%%%%%%%%%%%%%%%%%%%%%%%%%%%%%%%%%%%%%%%%

%%%%%%%%%%%%%%%%%%%%%%%%%%%%%%%%%%%%%%%%%%%%%%%%%%%%%%%%%%%%%%%%%%%%%%%%%%%%%%%%%%%%%%%%%%%%%
%\input{ch_numerics}

%%%%%%%%%%%%%%%%%%%%%%%%%%%%%%%%%%%%%%%%%%%%%%%%%%%%%%%%%%%%%%%%

\section{Conclusion}
We have studied the synchronization
of the velocities of different components in
a mixture flow
that is governed by 
a coupled system of isothermal Euler equations.
The coupling terms  appear as additional source terms that favor  the velocity 
alignment. 
An important application of the model is hydrogen blending in the natural gas flow
throught gas pipelines.

We have shown that under suitable  regularity conditions for the solution
and for appropriate boundary conditions  the 
$L^2$-norm of the difference of
the velocities  decays exponentially fast
as long as it is above a certain threshold. 
The decay rate is 
%up to a certain level
%that is 
governed by the order of the coupling terms.  
The  threshold for the synchronization error
depends on the reciprocal value of the order of the coupling  terms 
and on the size of the spatial pressure variations,
that is the partial derivatives of the pressures.
The proof is based on a  Lyapunov functions 
that has been chosen motivated by relative energy techniques.
Our result can be seen as a
theoretical justification
for the drift-flux models for
mixture flow. 
An extension of the result  to the case 
of networked systems is of interest
for future research. 
For the case of  networked systems,
the choice of suitable coupling conditions
at the junctions is crucial, see
for example \cite{rinaldo}.

{%\color{blue}
In our  analysis we have worked with classical solutions
of the dynamics that is governed by 
a system of quasi-linear balance laws.
Solutions of this type can break down after finite time,
therefore it is of high interest to 
investigate whether the
synchronization of the velocities still occurs
for weak entropy solutions that are less regular.
In \cite{wellposedness}, well-posedness 
results for networked systems of balance laws are given
for the case $n=1$, i.e. no mixture.
We expect that these results can be generalized to the case of the mixture flow.
The Lyapunov function has to be changed to be suitable for less regular solutions.
An obvious choice would be
\begin{equation}
   % \label{lerepresentationmix}
{\cal G}^e(t) = %\frac{1}{2}  \, 
\int_0^{L} 
 \sum_{i=1}^n 
%\frac{1}{2} \,  
\rho^i \, |  v^i -  v| \, dx
=
\int_0^{L} 
 \sum_{i=1}^n 
%\frac{1}{2} \,  
 | \rho^i \, v^i - \rho^i \, v| \, dx.
\end{equation}
However, this investigation is out of the scope of  this paper.
}

\textbf{Acknowledgements:}
This work was
supported by DFG in the 
%framework of the 
 Collaborative Research Centre
CRC/Transregio 154,
Mathematical Modelling, Simulation and Optimization Using the Example of Gas Networks,
Project C03 and  Project  C05,  Projektnummer 239904186.
%The authors thank the Bundesministerium für Bildung und Forschung (BMBF)
%for support under DAAD grant 57654073 'Uncertain data in control of PDE systems'.

%\bibliographystyle{siam}
\bibliographystyle{siamplain}

\end{document}